\documentclass[12pt]{amsart}

\usepackage[margin=1in]{geometry}  
\usepackage{graphicx}              
\usepackage{amsmath}               
\usepackage{amsfonts}              
\usepackage{amsthm}                
\usepackage{esint}
\usepackage{xcolor}
\usepackage{hyperref}		

\newtheorem{thm}{Theorem}[section]
\newtheorem{lem}[thm]{Lemma}
\newtheorem{prop}[thm]{Proposition}
\newtheorem{cor}[thm]{Corollary}

\newtheorem{rmk}[thm]{Remark}
\newtheorem{definition}[thm]{Definition}
\newtheorem{question}[thm]{Question}

\renewenvironment{abstract}{%
  \noindent\bfseries\abstractname:\normalfont}{}

\newcommand{\RR}{\mathbb{R}}      

\begin{document}

\title{Unique Continuation on Convex Domains}

\author{Sean McCurdy\\}

\maketitle

\begin{abstract}
In this paper, we obtain estimates on the quantitative strata of the critical set of non-trivial harmonic functions $u$ which vanish continuously on $V \subset \partial \Omega$, a relatively open subset of the boundary of a convex domain $\Omega \subset \mathbb{R}^n$.  In particular, these estimates improve dimensional estimates on $\{|\nabla u| =0\}$ both in $V \subset \partial \Omega$  and as it \textit{approaches} $V \cap \overline{\Omega}.$  These estimates are not obtainable by naively combining interior and boundary estimates and represent a significant improvement upon existing results for boundary analytic continuation in the convex case.
\end{abstract}

\tableofcontents

\section{Introduction}

Unique continuation is a fundamental property for functions which solve the Laplace and related linear equations.  A closely related problem is that of boundary unique continuation: given a domain $\Omega \subset \mathbb{R}^n$ and a function $u$ which is harmonic in $\Omega$ and vanishes continuously on $V \subset \partial \Omega$, how large can the set $\{Q \in V: |\nabla u| = 0\}$ be if $u \not \equiv 0$? Boundary unique continuation is closely tied to the Cauchy problem and questions of well-posedness and stability of solutions to boundary value problems (see, for instance, \cite{Tataru02} and \cite{AlessandriniRondiRossetVessella09}). In this paper, we address two questions.  First, we address the question of boundary unique continuation for harmonic functions on convex domains. Second, we also address the question of how the critical set $\{|\nabla u|=0\} \cap \Omega$ approaches $V \subset \partial \Omega$. 
We follow the approach of Garofalo and Lin \cite{GarofaloLin87} insofar as we make essential use of the Almgren frequency function.
And, because we want to obtain results on the full critical set $\{|\nabla u|=0\}$, we use packing estimates inspired by Cheeger, Naber, and Valtorta \cite{CheegerNaberValtorta15}.  These tools allow us to obtain results the strata of the critical set $\{ p \in \Omega : |\nabla u(p)| = 0\}$ \textit{as it approaches} $V \subset \partial \Omega.$

\subsection{Background on Boundary Unique Continuation for Harmonic Functions}  For dimensions $n \ge 3$, Bourgain and Wolff \cite{BourgainWolff90} have constructed an example of a function, $u: \RR^{n}_+ \rightarrow \RR$, which is harmonic in $\RR^n_+$, $C^1$ up to the boundary $\RR^{n-1} \subset \RR^n$, and for which both $u$ and $\nabla u$ vanish on a set of positive surface measure.  This result has been generalized by Wang \cite{Wang95}
to $C^{1, \alpha}$ domains, $\Omega \subset \RR^n,$ for $n \ge 3$.  However, the sets of positive measure for which these functions vanish are \textit{not} open.  

In general, the following question posed by Lin in \cite{LinFH91} is still open.

\begin{question}\label{open question}
Let $n \ge 2$ and $\Omega \subset \RR^n$ be an open, connected Lipschitz domain.  If $u$ is a harmonic function which vanishes continuously on a relatively open set $V \subset \partial \Omega$, does
\begin{align*}
\mathcal{H}^{n-1}(\{x \in V: |\nabla u| = 0\}) > 0
\end{align*}
imply that $u$ is the zero function?
\end{question}

If $u$ is non-negative, the techniques of PDEs on non-tangentially accessible (NTA) domains give a comparison principle \cite{Dahlberg77} which allows us to say that the norm of the normal derivative is point-wise comparable to the density of the harmonic measure with respect to the surface measure $d\sigma$.  Additionally, for Lipschitz domains it is well-known that the harmonic measure is mutually absolutely continuous with respect to $d\sigma$.  These two facts then imply that if the normal derivative vanishes on a set of positive (surface) measure, then $u$ must be identically $0$.  

The challenge is for harmonic functions $u$ which change sign.  For such functions, the aforementioned techniques fail completely because we cannot apply the Harnack principle.  Authors have therefore approached this problem by asking for additional regularity.  In \cite{LinFH91}, Lin proves that for $C^{1, 1}$ domains, $\Omega \subset \RR^n$, for $n \ge 2$, if $u$ is a non-constant harmonic function which vanishes on an open set $V \cap \partial \Omega$, then $\dim_{\mathcal{H}}(\{x \in V \cap \partial \Omega: |\nabla u| = 0\}) \le n-2.$  Similar results were later shown by Adolfsson and Escauriaza for domains with locally $C^{1, \alpha}$ boundary \cite{AdolfssonEscauriaza97}.  Relatedly, Kukavica and Nyst\"om showed that $\mathcal{H}^{n-1}(\{x \in V: |\nabla u| = 0\}) > 0$ implies that $u \equiv 0$ if $\partial \Omega$ is $C^{1}$ Dini \cite{KukavicaNystrom98}.  Recently, this result has been greatly improved. Kenig and Zhou \cite{KenigZhou21} employed powerful technicues from \cite{NaberValtorta17-1} and \cite{deLellisMarcheseSpadaroValtorta16}, have shown that for $C^1$ Dini domains, the $(n-2)$-generalized singular set $\{u = 0 = |\nabla u|\}$ has finite $(n-2)$-dimensional upper Minkowski content.

For merely convex domains $\Omega \subset \RR^n$, Adolfsson, Escauriaza, and Kenig showed that if $u$ is a harmonic function in $\Omega$ which vanishes continuously on a relatively open set $V \subset \partial \Omega$, then if $\{x \in V \cap \partial \Omega: |\nabla u| = 0\}$ has positive surface measure, $u$ must be a constant function \cite{AdolfssonEscauriazaKenig95}.  The method of attack pursued in \cite{AdolfssonEscauriazaKenig95} (and \cite{LinFH91}, \cite{AdolfssonEscauriaza97}, \cite{KukavicaNystrom98}) was centered on showing that the harmonic function is ``doubling" on the boundary in the following sense.  If $\Omega \subset \RR^n$, then there exists an absolute constant $M<\infty$ such that for all $B_{2r}(Q) \cap \partial \Omega \subset V$ 
\begin{align*}
\int_{B_{2r}(Q) \cap \Omega}u^2 dx \le M \int_{B_r(Q)\cap \Omega} u^2 dx.
\end{align*}

This doubling property allows the authors to show that the normal derivative is an $A_2$-Muckenhoupt weight with respect to surface measure, a kind of quantified version of mutual absolute continuity.  It is well-known that if $u$ vanishes in a surface ball $\Delta_r(Q)$ and the normal derivative of $u$ is a $2$-weight with respect to surface measure, then either $\{Q' \in \Delta_r(Q) : | \nabla u| = 0 \}$ has measure zero or $\{Q' \in \Delta_r(Q) : |\nabla u| > 0 \}$ has measure zero.  The improvement from measure to dimension bounds in \cite{AdolfssonEscauriaza97} and \cite{LinFH91} comes from applying an additional Federer dimension-reduction type argument.  

Recently, Tolsa has answered Question \ref{open question} in the affirmative for $C^1$ domains and Lipschitz domains with small Lipschitz constant \cite{Tolsa20}.

In this paper, we restrict our investigation to convex domains $\Omega$ and obtain bounds upon the \emph{full} generalized critical set $\{|\nabla u| = 0\}$. We note that this includes $\{x \in \partial \Omega: |\nabla u| = 0\}$ and $\{x \in \Omega: |\nabla u|=0\}$.  

\begin{thm}\label{e-reg}
Let $\Omega$ be a convex domain and $u \in C^{0}(\overline{\Omega})$ be a non-constant function such that $\Delta u = 0$ in $\Omega$.  Let $V \subset \partial \Omega$ be a relatively open set. If $u = 0$ on $V$, then for any compact subset $K \subset V$, there exists a radius $0<r(K)$ such that
\begin{align*}
& \dim_{\mathcal{H}}\left(B_{r}(K) \cap \{ |\nabla u| = 0\} \setminus \emph{sing}(\partial \Omega) \right)\\
& \le \overline{\dim_{\mathcal{M}}}\left(B_{r}(K) \cap \{ |\nabla u| = 0\} \setminus \emph{sing}(\partial \Omega)\right) \le n-2.
\end{align*}
Furthermore
\begin{align*}
\overline{\dim_{\mathcal{M}}}\left(B_{r}(K) \cap \overline{\{ |\nabla u| = 0\} \cap \Omega}\right) \le n-2.
\end{align*}
\end{thm}

The content of Theorem \ref{e-reg} is two-fold.  First, consider the results restricted to the boundary $\{|\nabla u|  = 0\} \cap V \subset \partial \Omega$.  Alberti \cite{Alberti94} proved (among other things) that the singular set of a convex function is a $C^2$ $(n-2)$-rectifiable set, which implies that the geometric singular set of a convex body satisfies $\dim_{\mathcal{H}}(\text{sing}(\partial \Omega)) \le n-2$.  Thus, Theorem \ref{e-reg} combined with \cite{Alberti94} implies 
\begin{align*}
\dim_{\mathcal{H}}\left(B_{r}(K) \cap \{ |\nabla u| = 0\}\right) \le n-2,
\end{align*}
which gives a strong improvement on the results of \cite{AdolfssonEscauriazaKenig95}, which proved that in this situation $\mathcal{H}^{n-1}(V \cap \{|\nabla u|=0\})=0.$

Second, Theorem \ref{e-reg} provides new insight into how the critical set interacts with $\partial \Omega.$  Returning to \cite{Alberti94}, Alberti also proved that the singular set of a convex function may be prescribed to be \emph{any} $C^2$ $(n-2)$-rectifiable set.  Thus, it can happen that $\overline{\dim}_{\mathcal{M}}(\text{sing}(\partial \Omega))=n-1$. On the other hand, from the interior perspective, \cite{NaberValtorta15} proved finite $(n-2)$-dimensional upper Minkowski content bounds on $\{|\nabla u|=0\}$ in the interior. But, naive application of these estimates degenerate as one approaches the boundary because upper Minkowski dimension is \textit{not} stable under countable unions.  Considering $\{|\nabla u| = 0 \} \cap \Omega$ as it approaches $\partial \Omega$, it was unknown whether or not $\{|\nabla u| = 0 \} \cap \Omega$ could oscillate wildly and have positive $(n-1)$-upper Minkowski dimension like $\text{sing}(\partial \Omega)$, would remain $(n-2)$-upper Minkowski dimensional like the interior, or if something in between these two held.  Theorem \ref{e-reg} proves that the set $\{|\nabla v|=0\} \cap \Omega$ cannot oscillate too wildly as it approaches $\partial \Omega$ and $\overline{\{|\nabla u|=0\} \cap \Omega} \cap K$ inherits its upper Minkowski dimension bounds from the interior rather than the boundary.

It is still an open question whether or not the $(n-2)$-upper Minkowski content of $\{|\nabla u|=0\} \cap \partial \Omega \setminus \text{sing}(\partial \Omega)$ is finite.

The author would like to thank Tatiana Toro, whose direction, advice, patience, and support can only be described as \textit{sine qua non}.  Additional thanks are due to Zihui Zhou for gently pointing out the errors in a previous version of this project and for a very thorough and patient reviewer, whose comments have greatly improved the presentation. 

\section{Definitions and Main Results}\label{S:defs}

Theorem \ref{e-reg} is a corollary to the Technical Theorem (Theorem \ref{T: main theorem 1}) and a containment result (Lemma \ref{e-reg containment}).  In order to state these results, we need the following definitions.

We define the following class of domains.
\begin{definition}(A normalized class of convex domains) \label{domain def}
 Let
$\mathcal{D}(n)$ be the collection of connected, open domains $\Omega \subset \RR^n$ which satisfy the following conditions:
\begin{enumerate}
\item $0 \in \partial \Omega$.
\item $\Omega \cap B_2(0)$ is convex.
\item $\Omega \cap (B_2(0))^c \not = \emptyset$.
\end{enumerate}
\end{definition} 

One of the key tools of this paper will be an Almgren frequency function, introduced by Almgren in \cite{Almgren79}.

\begin{definition}(Almgren frequency function) \label{N def}  Let $r >0$, $\Omega \subset \RR^n$, $u: \RR^n \rightarrow \RR$, such that $u \in C(B_{2r}(p)) \cap W^{1, 2}(B_{2r}(p))$and $p \in \overline \Omega$. We define the following quantities:
\begin{align*}
H_{\Omega}(p, r, u) &:= \int_{\partial B_r(p) \cap \overline{\Omega}} |u - u(p)|^2 d\sigma.\\
D_{\Omega}(p, r, u) &:= \int_{B_r(p)} |\nabla u|^2dx.\\
N_{\Omega}(p, r, u) &:= r \frac{D_{\Omega}(p, r, u)}{H_{\Omega}(p, r, u)}.
\end{align*}
\end{definition}

\begin{rmk}(Invariances of the Almgren frequency function) \label{R: rescalings}
This normalized version of the Almgren frequency function is invariant in the following senses.  Let $a, b, c \in \RR$ with $a, r \not = 0$.  If $w(x) = au(bx + p) + c$ and $T_{p,b}\Omega = \frac{1}{b}(\Omega - p)$ then 
$$
N_{\Omega}(p,r, u) = N_{T_{p, b}\Omega}(0, b^{-1}r, w).
$$
\end{rmk}

We now define the class of functions in which we will work in this paper.
\begin{definition}(A class of functions) \label{A def}
 Let
$\mathcal{A}(n, \Lambda)$ be the set of functions, $u: \RR^n \rightarrow \RR$, which have the following properties:
\begin{enumerate}
\item $u: \RR^n \rightarrow \RR$ is harmonic in a convex domain, $\Omega \in \mathcal{D}(n)$.
\item $u \in C(\overline{B_2(0)})$, $u = 0$ on $\Omega^c \cap B_2(0)$, and $u$ is non-constant.
\item $N_{\Omega}(0, 2, u) \le \Lambda$.
\end{enumerate}
\end{definition} 

We shall use rescalings which are adapted to the quantitative stratification methods introduced by Cheeger and Naber in \cite{CheegerNaber13} for studying the regularity of stationary harmonic maps and minimal currents.  

\begin{definition}(Rescalings) \label{L2 rescaling def}
 Let $u \in \mathcal{A}(n, \Lambda)$, and let $\Omega$ be its associated domiain.  We define the rescaled function, $T_{x, r}u$ of $u$ at a point $x \in B_{1}(0)$ at scale $0<r<1$ by
 \begin{align*}
     T_{x, r}u(y) := \frac{u(x +ry) - u(x)}{\left(\int_{\partial B_1(0)\cap \overline{\Omega}} (u(x + ry) -u(x))^2d\sigma(y)\right)^{1/2}}.
 \end{align*}
In the case that the denominator is zero, we define $T_{x, r}u= \infty$.

We shall break with established convention and denote the rescalings of $\Omega$ in analogy with the rescaling of functions. Let $T_{p, r}\Omega := \frac{\Omega - p}{r}$ and $T_{p, r}\partial \Omega := \frac{\partial \Omega - p}{r}.$
\end{definition}

The geometry we wish to capture with the rescalings $T_{x, r}f$ are encoded in their translational symmetries.
 
\begin{definition}(The class of blow-up profiles)\label{symmetric def}
Let $u \in C(\RR^n)$. We say $u$ is $0$-\emph{symmetric} if $u$ satisfies one of the following conditions.
\begin{enumerate}
\item $u$ is a homogeneous harmonic polynomial. 
\item $u(x) = \phi(x-p) + c$ for some function $\phi$ which is homogeneous and harmonic in a convex cone, $\Omega' \in \mathcal{D}(n)$, some point $p \in \mathbb{R}^n \cap \Omega'$, and some $c \in \mathbb{R}$.
\end{enumerate}
We will say that $u$ is $k$-\emph{symmetric} if $u$ is $0$-symmetric and there exists a $k$-dimensional subspace $V$ such that $u(x+y) =u(x)$ for all $x \in \RR^n$ and all $y \in V$.
\end{definition}

We now define the quantitative version of symmetry which describes how close to being $k$-symmetric a function is in a ball, $B_r(x) \subset \mathbb{R}^n$.

\begin{definition}(Quantitative symmetry) \label{quant symmetric def}
 For any $u \in \mathcal{A}(n, \Lambda)$ with associated domain $\Omega$, $u$ will be called $(k, \epsilon, r, p)$-\emph{symmetric} if there exists a $k$-symmetric function $P$ such that 
 \begin{itemize}
 \item[1.] $\int_{\partial B_1(0)} |P|^2d\sigma = 1$.
 \item[2.] $\int_{B_1(0) \cap \overline{T_{p, r}\Omega}}|T_{p, r}u - P|^2 dV < \epsilon.$
 \end{itemize}
 
We shall say that $u$ is $(k, \delta_0, r, p)$-symmetric \emph{with respect to a $k$-dimensional subspace $V$} if there is a $k$-symmetric function $P$ which verifies that $u$ is $(k, \delta_0, r, p)$-symmetric such that $P(x+y) = P(x)$ for all $y \in V$.
\end{definition}

 \begin{definition}(Quantitative Generalized Critical Strata)\label{quant strat def}
 Let $u \in \mathcal{A}(n, \Lambda)$ with $\Omega \in \mathcal{D}(n)$ its associated domain.  For $0<\epsilon$, $0< r \le 1$, and integer $0 \le k \le n-1$ we denote the $(k, \epsilon, r)$-\emph{generalized critical strata} of $u$ by $\mathcal{C}^k_{\epsilon, r}(u)$, and we define it by 
\begin{equation*}
\mathcal{C}^k_{\epsilon, r}(u) := \{x \in \overline{\Omega} : u \text{  is not  } (k+1, \epsilon, s, x) \text{-symmetric for all  } r \le s \le 1\}.
\end{equation*}
We shall also use the notation $\mathcal{C}^k_{\epsilon}(u)$ for $\mathcal{C}^k_{\epsilon, 0}(u).$
\end{definition}

The quantitative strata of the generalized critical set behave very well under $L^2$-convergence of functions.  That is, if $u_i \in \mathcal{A}(n, \Lambda)$ and $x_i \in B_1(0) \cap \overline{\Omega_i}$ such that $x_i \rightarrow x \in \overline{B_1(0)}$, then if $u_i$ is $(k, \epsilon, r_i, x_i)$-symmetric and $T_{x_i, 1}u_i \rightarrow u_{\infty}$ in $L^2(B_2(0))$ then $u_{\infty}$ is $(k, \epsilon, 1, 0)$-symmetric. As a special case, this implies that the $C^{k}_{\epsilon, r}(u)$ are closed.

It is immediate from the definitions that $\mathcal{C}^k_{\epsilon, r}(u) \subset \mathcal{C}^{k'}_{\epsilon', r'}(u)$ if $k \le k', \epsilon' \le \epsilon, r \le r'$.  

\begin{definition}(Qualitative Generalized Critical Set)
Let $u \in \mathcal{A}(n, \Lambda)$ with $\Omega \in \mathcal{D}(n)$ its associated domain. Using the quantitative generalized critical strata, we define the \emph{generalized critical set} of $u$    $\mathcal{C}^{n-2}(u) := \cup_{\eta} \cap_{r} \mathcal{C}^{n-2}_{\eta, r}(u)$. In turn, we define the strata of the generalized critical set as follows
$\mathcal{C}^k(u) := \cup_{\eta} \cap_{r} \mathcal{C}^k_{\eta, r}(u).$
\end{definition}

We shall use the convention that for any $A \subset \RR^n$,  $B_r(A) = \{x \in \RR^n : d(A, x) < r \}$.  Recall that we can define upper Minkowski $s$-content by
\begin{align}\label{upper minkowski content}
\mathcal{M}^{*,s}(A) = \limsup_{r \rightarrow 0}\frac{\text{Vol}(B_r(A))}{\omega_{n-s}r^{n-s}}
\end{align}
and upper Minkowski dimension as 
$$\overline \dim_{\mathcal{M}}(A) = \inf\{s: \mathcal{M}^{*s}(A) = 0 \} = \sup \{s: \mathcal{M}^{*s}(A) > 0 \}.$$

Now, we can state the main technical results.
\begin{thm}(Technical Theorem)\label{T: main theorem 1}
Let $u \in \mathcal{A}(n, \Lambda)$, then for any $r_0> 0$ and all $0< r_0< r$
\begin{equation}\label{e:main theorem 1 estimate}
Vol\left(B_r\left(\mathcal{C}^k_{\epsilon, r_0}(u) \cap B_{\frac{1}{4}}(0)\right)\right) \le C(n, \Lambda, k, \epsilon)r^{n-k-\epsilon}.
\end{equation}
In particular, letting $r_0 \rightarrow 0$
\begin{equation*}
\mathcal{H}^{k+\epsilon}\left(\mathcal{C}^{k}_{\epsilon}(u) \cap B_{\frac{1}{4}}(0)\right) \le \mathcal{M}^{*, k+\epsilon}\left(\mathcal{C}^{k+\epsilon}_{\epsilon}(u) \cap B_{\frac{1}{4}}(0)\right) \le C(n, \Lambda, k, \epsilon)
\end{equation*}
\end{thm}

\begin{lem}(Containment)\label{e-reg containment}
There exists an $0<\epsilon = \epsilon(n, \Lambda)$ such that
\begin{align}\label{e:containment 1}
\left(\overline{\mathcal{C}^{n-2}(u) \cap \Omega}\right) \cap B_{\frac{1}{8}}(0)\subset \mathcal{C}^{n-2}_{\epsilon}(u)
\end{align}
and
\begin{align}\label{e:containment 2}
\left(\mathcal{C}^{n-2}(u) \cap \partial \Omega \setminus \emph{sing}(\partial \Omega) \right) \cap B_{\frac{1}{8}}(0)\subset \mathcal{C}^{n-2}_{\epsilon}(u)
\end{align}
\end{lem}

\textit{Proof of Theorem \ref{e-reg} assuming Lemma \ref{e-reg containment}}  For each point $x \in K \subset V,$ there is a radius $0<r$ such that $B_{4r}(x) \cap \partial \Omega \subset V.$  Since $K$ is compact, we may find a finite subcover $\{B_{r_i}(x_i)\}_i.$  Thus, Lemma \ref{e-reg containment} implies that in each $B_{r_i}(x_i)$
\begin{align*}
    \inf\{s: \mathcal{M}^{*,s}(B_{r_i}(x_i) \cap \{|\nabla u| = 0\} \setminus \text{sing}(\partial \Omega)) < \infty\} \le n-2.
\end{align*}
Since upper Minkowski dimension is stable under finite unions, the first claim of Theorem \ref{e-reg} holds. The second follows from an identical argument using Lemma \ref{e-reg containment}. \qed

\subsection{Outline of the Paper.} 
The structure of this paper is roughly in four parts.  Section \ref{S:monotonicity} and Section \ref{sec:uniform frequency bound Singular set} use the geometric techniques of \cite{GarofaloLin87}, \cite{HanLin_nodalsets} (and many, many others) to establish that the Almgren frequency is monotonically non-decreasing and bounded on $\{u=0\}$.  Section \ref{sec: compactness Zero Set} uses these results to establish compactness of $\{T_{p, r}u\}$ for $u(p) = 0$.  Section \ref{S: off zero set} extends these results to $p \in \Omega$ such that $u(p)\not= 0$.

The second part of this paper is devoted to obtaining geometric control upon $\mathcal{C}^k_{\epsilon, r}(u).$  The general idea is to employ the usual ``frequency pinching" (Lemma \ref{quant rigidity}) and cone-splitting  results (Lemma \ref{cone splitting}). However, because we are considering $N_{\Omega}(p, r, u)$ at points $p$ such that $u(p) \not = 0$, the Almgren frequency is \textit{not} monotonic. This is overcome by proving that if $\text{dist}(p, \{u=0 \}) << r$, then \begin{align*}
N_{\Omega}(p, 1, u)- N_{\Omega}(p, 1/10, u) \le \gamma, 
\end{align*}
for $0< \gamma$ sufficiently small implies that $u$ is $(0, \epsilon, 1, p)$-symmetric.  In Corollary \ref{c: application of cone splitting}, we prove that if $x \in \mathcal{C}^{k}_{\epsilon}(u)$ then $\{p \in \overline{\Omega}: u \text{ is } (0, \delta, r, p)\text{-symmetric}\} \cap B_r(x)$ is contained in a tubular neighborhood of a $k$-plane $L^k$.

The third part of this paper is devoted to obtaining packing estimates to prove Theorem \ref{T: main theorem 1}. To do so, we use the tools of \cite{CheegerNaberValtorta15}, which do not require restricting to a level set or the delicate machinery which powers the finer estimates of \cite{deLellisMarcheseSpadaroValtorta16}.  The fact that we do not control the tilt of approximating $L^k$ at different scales accounts for the $(k+\epsilon)$-dimensional results.

The fourth part of this paper is devoted to proving the containment results which prove Lemma \ref{e-reg containment}.  

Throughout this paper, the constant $C$ will by ubiquitous and represent different constants even within the same string of inequalities.  A constant written $C(n, \Lambda)$ will only depend upon $n$ and $\Lambda$, but each instantiation may represent a distinct constant.

\section{The Almgren frequency function}\label{S:monotonicity}

In this section, we develop crucial properties of the Almgren frequency function.  The main results of this Section are the monotonicity of the Almgren frequency on $p \in \{u=0\} \cap \overline{\Omega}$.

We now note some of the elementary properties of $H_{\Omega}(p, r, u), D_{\Omega}(p, r, u), N_{\Omega}(p, r, u),$ and their derivatives.

\begin{lem}
Let $u \in \mathcal{A}(n, \Lambda)$ and $p \in \overline{\Omega} \cap B_1(0)$ and all $0< r< 1$ 
\begin{equation}\label{H derivative}
\frac{d}{dr}H_{\Omega}(p, r, u) = \frac{n-1}{r}H_{\Omega}(p, r, u) + 2D_{\Omega}(p, r, u) + 2u(p) \int_{\partial \Omega \cap B_r(p)} \nabla u \cdot \vec \eta d\sigma
\end{equation}
\begin{align} \nonumber \label{D derivative}
\frac{d}{dr}D_{\Omega}(p, r, u) = & \frac{n-2}{r}D_{\Omega}(p, r, u) + 2\int_{\partial B_r(p)}(\nabla u \cdot \vec \eta)^2d\sigma\\
&  + \int_{\partial \Omega \cap B_r(p)}(Q-p) \cdot \vec \eta (\nabla u \cdot \vec \eta)^2 d\sigma(Q)
\end{align}
\begin{align} \label{L: avg log H derivative}
\frac{d}{dr}\ln(\frac{1}{r^{n-1}}H_{\Omega}(p, r, u)) & = \frac{2}{r}N_{\Omega}(p, r, u) +  2\frac{u(p) \int_{\partial \Omega \cap B_r(p)} \nabla u \cdot \vec \eta d\sigma}{H_{\Omega}(p, r, u)}\\
\frac{d}{dr}\ln(H_{\Omega}(p, r, u)) & = \frac{n-1}{r} + \frac{2}{r}N_{\Omega}(p, r, u) + 2 \frac{u(p) \int_{\partial \Omega \cap B_r(p)} \nabla u \cdot \vec \eta d\sigma}{H_{\Omega}(p, r, u)},
\end{align}
where $\vec \eta$ is the unit \textit{outer} normal of the relevant domains.
\end{lem}

\begin{proof} In the interior setting, for $\overline{B_r(p)} \subset \Omega$, these identities follow from straightforward computation.  (\ref{H derivative}) follows from the change of variables, $y \rightarrow rx + p$, and the divergence theorem.  (\ref{D derivative}) relies upon the Rellich-Necas Identity, 
\begin{equation}
\text{div}(X|\nabla u|^2) = 2\text{div}((X \cdot \nabla u)\nabla u) + (n-2)|\nabla u|^2,
\end{equation}
the divergence theorem, and the fact that $u$ vanishes on the boundary.  The last two equations follow immediately from \ref{H derivative}.  Without exception, the standard interior computations go through identically for radii for which $B_r(p) \cap \partial \Omega \not = \emptyset$, where we also use the identity 
\begin{align*}
\int_{B_r(p)}(u - u(p))\Delta u = u(p) \int_{\partial \Omega \cap B_r(p)} \nabla u \cdot \vec \eta d\sigma,
\end{align*}
where $\vec \eta$ is the \textit{outer} unit normal vector to $\Omega$. See \cite{HanLin_nodalsets} Theorem 2.2.3, Corollary 2.2.5 and \cite{AdolfssonEscauriazaKenig95} (Proof of the Doubling Property) for details.  
\end{proof}

The following lemma records a useful identity which follows from the previous lemma by straightforward computation.

\begin{lem}\label{N derivative calculation}
For $u \in \mathcal{A}(n, \Lambda)$ and $p \in \overline{\Omega} \cap B_1(0)$ and all $0< r< 1$, $\frac{d}{dr}N_{\Omega}(p, r, u)$ may be decomposed into four terms
\begin{align}
\frac{d}{dr}N_{\Omega}(p, r, u) = & N_1'(r) + N_2'(r) + N_3'(r) + N_4'(r) ,
\end{align}
where 
\begin{align*}
N_1'(r) := & \frac{1}{H_{\Omega}(p, r, u)^2} 2r\left[ H_{\Omega}(p, r, u) \int_{\partial B_r(p) \cap \overline{\Omega}}(\nabla u \cdot \vec \eta)^2d\sigma - \left(\int_{\partial B_r(p) \cap \overline{\Omega}}(u - u(p))(\nabla u \cdot \vec \eta)d\sigma\right)^2\right]\\
N_2'(r) := & \frac{1}{H_{\Omega}(p, r, u)} r \int_{\partial \Omega \cap B_r(p)}(Q-p) \cdot \vec \eta (\nabla u \cdot \vec \eta)^2 d\sigma(Q)  \\
N_3'(r) := & 2N_{\Omega}(p, r, u) \frac{u(p)}{H_{\Omega}(p, r, u)} \int_{\partial \Omega \cap B_r(p)} \nabla u \cdot \vec \eta d\sigma\\
N_4'(r) := & \frac{2r}{H_{\Omega}(p, r, u)^2}  \left(u(p) \int_{\partial \Omega \cap B_r(p)} \nabla u \cdot \vec \eta d\sigma\right)^2,
\end{align*}
and $\vec \eta$ is the unit \textit{outer} normal.
\end{lem}

\begin{lem}\label{N_1 redefine}
Let $u \in \mathcal{A}(n, \Lambda)$ and $p \in \overline{\Omega} \cap B_1(0) \cap \{u=0\}$ and all $0< r< 1$,
\begin{align*}
N_1'(r) = & \frac{2r\left[ H_{\Omega}(p, r, u) \int_{\partial B_r(p)}(\nabla u \cdot \vec \eta)^2d\sigma - \left(\int_{\partial B_r(p)}(u - u(p))(\nabla u \cdot \vec \eta)d\sigma\right)^2\right]}{H_{\Omega}(p, r, u)^2}\\
 = &  \frac{2}{r H_{\Omega}(p, r, u)} \left( \int_{\partial B_r(p) \cap \overline{\Omega}} |\nabla u \cdot (y-p) - N_{\Omega}(p, r, u) (u(y)-u(p))|^2d\sigma(y)\right).
\end{align*}
\end{lem}

\begin{proof}
Recall that by the Cauchy-Schwarz inequality, we have that for $\lambda = \frac{\langle w, v \rangle}{||v||^2}$
$$
||v ||^2 ||w-\lambda v ||^2 = |w|||^2 ||v||^2 - |\langle w, v \rangle|^2.
$$
 
Choosing $w = \nabla u \cdot (y-p)$ and $v = u- u(p)$, we have
\begin{align*}
N_1'(r) = & H_{\Omega}(p, r, u)^{-1}2r \left( \int_{\partial B_r(p)\cap \overline{\Omega}} |(u)_{\nu} - \frac{1}{r}\lambda(p, r, u) (u-u(p))|^2d\sigma \right) \\
 = &  \frac{2}{r H_{\Omega}(p, r, u)} \left( \int_{\partial B_r(p)\cap \overline{\Omega}} |\nabla u \cdot (y-p) - \lambda(p, r, u) (u(y)-u(p))|^2d\sigma(y)\right),
\end{align*} 
where 
\begin{align}\label{lambda def}
\lambda(p, r, u) := \frac{\int_{\partial B_{r}(p)\cap \overline{\Omega}} (u(y) - u(p))\nabla u \cdot (y-p) d\sigma(y)}{H_{\Omega}(p, r, u)}.
\end{align}
The Divergence theorem then implies that for $u(p)=0$, $\lambda(p, r, u) = N_{\Omega}(p, r, u).$
\end{proof}

\begin{lem}(Monotonicity)\label{N monotonicity 1}
Let $u \in \mathcal{A}(n, \Lambda)$ and $p \in \{u = 0\} \cap \overline{\Omega} \cap B_1(0)$, then $N_{\Omega}(p, r, u)$ is monotonically non-decreasing in $0< r < 1$. 
\end{lem}

\begin{proof} 
Recall Lemma \ref{N derivative calculation}.  Note that $N_1'(r)$ is non-negative by the Cauchy-Schwartz inequality.  Furthermore, because $N_2'(r)$ is non-negative because $\partial \Omega$ is a convex surface and $\vec \eta$ is the outer normal, $(Q- p) \cdot \vec \eta \ge 0$ for all $p \in \Omega \cap B_2(0)$ and all $Q \in \partial \Omega \cap B_2(0).$ Observe that $N_3'(r) = N_4'(r) = 0$ because $u(p) = 0$.  Therefore, $\frac{d}{dr}N_{\Omega}(p, r, u)$ is non-negative.
\end{proof}

\section{The Zero set: Uniform Frequency Bounds}\label{sec:uniform frequency bound Singular set}
The main result in this section is Lemma \ref{N bound lem 1}, which gives a uniform bound on the Almgren frequency function for all $p \in \overline{\Omega} \cap B_{\frac{1}{4}}(0)$ for which $u(p)=0$ and all $0< r \le \frac{1}{2}$.  We begin with a few basic results.

\begin{lem} \label{L: H doubling-ish 1}($H_{\Omega}(p, r, u)$ is Doubling)
Let $u \in \mathcal{A}(n, \Lambda)$ with $p \in B_{1}(0) \cap \{u=0\} = \overline{\Omega}$.  For any $0 < s < S \le 1$
\begin{equation}
H_{\Omega}(p, S, u) \le \left(\frac{S}{s}\right)^{(n-1) + 2N_{\Omega}(p, S, u)} H_{\Omega}(p, s, u).
\end{equation}
\end{lem}

\begin{proof} Recalling Equations (\ref{H derivative}) and (\ref{L: avg log H derivative})
\begin{align*}
\ln\left(\frac{H_{\Omega}(p, S, u)}{H_{\Omega}(p, s, u)}\right) & =  \ln\left(H_{\Omega}(p, S,u)\right) - \ln\left(H_{\Omega}(p, s, u)\right)\\
& =  \int_s^S \frac{H_{\Omega}'(p, r, u)}{H_{\Omega}(p, r, u)}dr\\
 & =  \int_s^S \frac{n-1}{r} + \frac{2}{r}N_{\Omega}(p, r, u).
\end{align*}
We bound $N_{\Omega}(p, r, u)$ by $N_{\Omega}(p, S, u)$ using Lemma \ref{N monotonicity 1}.  Plugging in these bounds, we have that for $r \in [s, S]$
\begin{align*}
\ln\left(\frac{H_{\Omega}(p, S, u)}{H_{\Omega}(p, s, u)}\right) & \le [(n-1) + 2N_{\Omega}(p, S, u)] \ln(r)|^{S}_s.
\end{align*}
Evaluating and exponentiating gives the desired result.
\end{proof}

 \begin{rmk}\label{R: H doubling-ish}
 Because $N_{\Omega}(p, r, u)$ is monotonic for $p \in B_{1}(0) \cap \{u=0\} \cap \overline{\Omega}$, we can also extract the inequality
 \begin{align*}
\ln\left(\frac{H_{\Omega}(p, S, u)}{H_{\Omega}(p, s, u)}\right) & \ge [(n-1) + 2N_{\Omega}(p, s, u)] \ln(r)|^{S}_s,
\end{align*}
which leads to 
\begin{equation}
H_{\Omega}(p, s, u) \le \left(\frac{s}{S}\right)^{(n-1) + 2N_{\Omega}(p, s, u)} H_{\Omega}(p, S, u).
\end{equation}

If $S = 1$ and $u = T_{p, r}u$, then we have that for all $1> s >0$
\begin{equation}
H_{\Omega}(0, s, T_{p, r}u) \le s^{(n-1) + 2N_{\Omega}(Q, 0, T_{p, r}u)}.
\end{equation} 
 \end{rmk}
 
We are now ready for the main result of this section.

\begin{lem}(Uniform bound on frequency)\label{N bound lem 1}
Let $u \in \mathcal{A}(n, \Lambda)$, as above. There is a constant, $C_1(n, \Lambda)$ such that for all $p \in \{u=0\} \cap \overline{\Omega} \cap B_{\frac{1}{4}}(0)$ and all $r \in (0, \frac{1}{2})$ 
\begin{equation}
 N_{\Omega}(p, r, u)  \le  C_1(n, \Lambda).
\end{equation}
\end{lem}

\begin{proof} Recall that $0 \in \partial \Omega$ and that the Almgren frequency function is invariant under rescalings.  Therefore, we normalize our function $u$ by the rescaling $v = T_{0,1}u$.

Therefore, applying Lemma \ref{L: H doubling-ish 1} to $Q = 0$, letting $r = cR$, and integrating both sides with respect to $R$ from $0$ to $S$, we have that for any $c \in (0, 1)$
\begin{align*}
\int_{B_{S}(0)}|v|^2dV  &\le \int_0^{S} \left(\frac{1}{c}\right)^{(n-1) + 2N_{\Omega}(0, R, v)} \int_{\partial B_{cR}(0)}|v|^2d\sigma dR\\
& \le \left(\frac{1}{c}\right)^{(n-1) + 2N_{\Omega}(0, S, v)} \int_0^{S}  \int_{\partial B_{cR}(0)}|v|^2d\sigma dR.
\end{align*}
Thus, letting $S = 1$ and $c= \frac{1}{16}$ and dividing by $\omega_n$ we obtain
\begin{equation}\label{N bound base ineq 1}
\fint_{B_1(0)}|v|^2dV \le 16^{2N_{\Omega}(0, 1, v)} \fint_{B_{\frac{1}{16}}(0)} |v|^2dV.
\end{equation}

Thus, for any $p \in \{v=0\} \cap \overline{\Omega} \cap B_{\frac{1}{4}}(0)$ by inclusion
\begin{align*}
\int_{B_{1}(0)}|v|^2dV  &\ge \int_{B_{3/4}(p)}|v|^2dV, \qquad \int_{B_{\frac{1}{16}}(0)} |v|^2dV \le \int_{B_{\frac{9}{16}}(p)}|v|^2dV.
\end{align*}
Therefore, substituting these bounds into (\ref{N bound base ineq 1})
\begin{equation}\label{N bound base ineq 2}
\fint_{B_{\frac{3}{4}}(p)}|v|^2dV \le 16^{2N_{\Omega}(0, 1, v)}\left(\fint_{B_{\frac{9}{16}}(p)} |v|^2dV\right).
\end{equation}

Now, we wish to bound $\fint_{B_{\frac{3}{4}}(p)}|v|^2dV$ from below and $\fint_{B_{\frac{9}{16}}(p)} |v|^2dV$ from above.  By (\ref{H derivative}), if $v(p) = 0$, $\frac{d}{dr} \int_{\partial B_r(p)} |v|^2d\sigma \ge 0$.  Thus, for all $p \in \{v=0\} \cap \overline{\Omega} \cap B_{\frac{1}{4}}(0)$, we bound
\begin{align*}
\int_{B_{\frac{3}{4}}(p)}|v|^2dV & \ge \int_{\frac{5}{8}}^{\frac{3}{4}} \int_{\partial B_r(p)} |v|^2 d\sigma dr \ge  c \int_{\partial B_{\frac{5}{8}}(p)} |v|^2 d\sigma.\\
\int_{B_{\frac{9}{16}}(p)} |v|^2dV & \le \int_0^{\frac{9}{16}} \int_{\partial B_{\frac{9}{16}}(p)}|v|^2 d\sigma dr \le c\int_{\partial B_{\frac{9}{16}} (p)} |v|^2d\sigma.
\end{align*}
Plugging the above bounds into (\ref{N bound base ineq 2}) and dividing, we obtain 
\begin{align} \label{ugly 1}
\frac{\fint_{\partial B_{\frac{5}{8}}(p)} |v|^2 d\sigma}{\fint_{\partial B_{\frac{9}{16}}(p)}|v|^2 d\sigma} \le & C(n) 16^{2N_{\Omega}(0, 1, v)}
\end{align}
for all $p \in \{v=0\} \cap \overline{\Omega} \cap B_{\frac{1}{4}}(0)$

Recalling Equation (\ref{L: avg log H derivative}) and Lemma \ref{N monotonicity 1}, we see
\begin{align*}
\ln\left(\fint_{\partial B_{\frac{5}{8}}(p)} v^2d\sigma \right) - \ln\left(\fint_{\partial B_{\frac{9}{16}}(p)} v^2d\sigma \right) & = \int_{\frac{9}{16}}^{\frac{5}{8}} \frac{d}{dr}\ln\left(\frac{1}{r^{n-1}}H_{\Omega}(p, r, v)\right) dr\\
& =  \int_{\frac{9}{16}}^{\frac{5}{8}} \frac{2}{r}N_{\Omega}(p, r, v)dr\\
& \ge  2N_{\Omega}(p, \frac{1}{2}, v)\left(\ln\left(\frac{5}{8}\right) - \ln\left(\frac{9}{16}\right)\right) \\
& \ge  2cN_{\Omega}(p, \frac{1}{2}, v).
\end{align*}
Thus, (\ref{ugly 1}) gives us that
\begin{align*}
 2c\left[N_{\Omega}(p, \frac{1}{2}, v)\right] \le & \ln\left(\frac{\fint_{\partial B_{\frac{5}{8}}(p)} |v|^2d\sigma}{\fint_{\partial B_{\frac{9}{16}}(p)} |v|^2d\sigma}\right)\\
  \le & \ln\left(C(n) (16)^{2N_{\Omega}(0, 1, v)}\right)\\
  = & 2N_{\Omega}(0, 1, v)\ln(16) + C(n)\\
  \le & 2\Lambda \ln(16) + C(n).
\end{align*}

Now, Lemma \ref{N monotonicity 1}, gives that for $\frac{1}{2} > s > 0$ $N_{\Omega}(p, 1/2, v) \ge  N_{\Omega}(p, s, v)$. Since $N_{\Omega}(p, r, v) = N_{\Omega}(p, r, u),$ we have the desired claim.
\end{proof}

\section{The Zero Set: Compactness}\label{sec: compactness Zero Set}

The uniform bounds on the Almgren frequency function allow us to prove compactness results on the collection of rescaling $\{T_{p, r}u\}$.  The main results of this Section are weak compactness (Lemma \ref{compactness 1}), the geometric non-degeneracy of the domains $\Omega$ (Corollary \ref{C: geometric non-degeneracy}).

We now state a sequence of preliminary corollaries to Lemma \ref{N bound lem 1}.  We shall denote the $C^{0, \gamma}(B_1(0))$-norm by
$$
||u||_{C^{0, \gamma}(B_1(0))} := ||u||_{C^0(B_1(0))} + \sup_{\substack{x, y \in B_1(0)\\ x \not=y}} \frac{|u(x) - u(y)|}{|x-y|^{\gamma}}.$$

\begin{lem}(Uniform Holder continuity)\label{L: unif Holder bound 1}
Let $u \in \mathcal{A}(n, \Lambda)$, $Q \in B_{\frac{1}{4}}(0) \cap \partial \Omega$, and $r \in (0,1/2]$.  Then
\begin{equation}
|| T_{Q, r}u ||_{C^{0, \gamma}(B_1(0))} \le C(n, \Lambda).
\end{equation}
\end{lem}

We defer the proof of this statement to the Appendix A.  The techniques are standard.

\begin{cor}(Non-degeneracy of domains)\label{C: geometric non-degeneracy}
Let $u \in \mathcal{A}(n, \Lambda)$ and $\Omega \in \mathcal{D}(n)$ it's associated convex domain.  There exists a constant, $0<c = c(\Lambda, n)$ such that for all $Q \in \partial \Omega \cap B_{\frac{1}{4}}(0)$ and $0< r \le \frac{1}{2}$, $\partial B_1(0) \cap T_{Q,r}\Omega$ is a relatively open convex surface with 
$$\mathcal{H}^{n-1}(\partial B_1(0) \cap T_{Q,r}\Omega) > c.$$
\end{cor}

\begin{proof} That $\partial B_r(Q) \cap \Omega$ is relatively open and relatively convex is immediate from the definition of $\Omega.$  By Lemma \ref{L: unif Holder bound 1} we see that $\max_{B_1(0)}|T_{Q, r}u(x)| \le C(n, \Lambda)$ and $H_{T_{Q,r}\Omega}(0, 1, T_{Q, r}u) = 1$.  Furthermore, we have that 
\begin{align*}
H_{T_{Q,r}\Omega}(0, 1, T_{Q, r}u) & \le \mathcal{H}^{n-1}(\partial B_1(0) \cap T_{Q,r}\Omega) C^2.
\end{align*}
Therefore, $\mathcal{H}^{n-1}(\partial B_1(0) \cap T_{Q,r}\Omega) \ge C^{-2} = c.$
\end{proof}

\begin{cor}\label{sub-linear bound}
For all $u \in \mathcal{A}(n, \Lambda)$, $Q_0 \in \partial \Omega \cap B_{1/4}(0),$ and $0 < r \le \frac{1}{2}.$ the following estimate holds.  Let $Q \in T_{Q_0, r}\partial \Omega \cap B_{\frac{1}{2}}(0)$ and let $L_{Q}$ be a supporting hyperplane to $Q \in T_{Q_0, r}\partial \Omega$.  Then for all $p \in \overline{T_{Q_0, r}\Omega} \cap B_{\frac{1}{4}}(Q)$ 
\begin{align*}
|T_{Q_0, r}u(p)| \le C(n, \Lambda) \emph{dist}(p, L_Q).
\end{align*}
In particular, $|\nabla T_{Q_0, r}u(Q)| \le C(n, \Lambda)$ for all $Q \in \partial T_{Q_0, r}\Omega \cap B_{\frac{1}{2}}(0)$.
\end{cor}

\begin{proof}
Let $\mathbb{H}_Q$ be the half-space with boundary $L_Q$ which contains $T_{Q_0, r}\Omega.$ Consider the Dirichlet problem
\begin{align*}
\Delta \phi = &  0  \qquad \qquad \text{  in  } \mathbb{H}_Q \cap B_{\frac{1}{2}}(Q),\\
 \phi = &  \begin{cases}
 C(n, \Lambda) &  \text{on  } \partial B_{\frac{1}{2}}(Q) \cap T_{Q_0, r}\Omega\\
 0 & \text{on  } \partial(B_1(Q) \cap \mathbb{H}_Q) \setminus (\partial B_{\frac{1}{2}}(Q) \cap T_{Q_0, r}\Omega),
 \end{cases}
\end{align*}
where we choose $C(n, \Lambda)$ to be the same constant in Lemma \ref{L: unif Holder bound 1} for which we have $\sup_{\partial B_1(0)} |T_{Q_0, r}u| \le C(n, \Lambda)$.  Note that for any $Q \in \partial \Omega$ $\mathbb{H}_Q \cap B_{1/2}(Q)$ is a Weiner regular domain and the boundary data is piecewise continuous, so a unique solution $\phi$ must exist.

By the maximum principle $T_{Q_0, r}u \le \phi$ in $T_{Q_0, r}\Omega \cap B_{\frac{1}{2}}(Q)$.  We now argue that $\phi$ is comparable to a linear function in $B_{\frac{1}{4}}(Q) \cap \mathbb{H}_Q$.  
 
Let $L$ be the affine linear function with $\{L = 0 \} = L_Q$ such that 
$$\max_{\partial B_{1/2}(Q)}L= \max_{\partial B_{1/2}(Q)}\phi = C(n, \Lambda).
$$
By \cite{JerisonKenig82} Theorem 5.1., there is a constant $C$ such that we have that for all $x \in B_{\frac{1}{4}}(Q) \cap \mathbb{H}_Q$
\begin{align*}
 \phi(x) \le C L(x),
\end{align*}
where $C$ depends only upon the geometry of $B_{\frac{1}{4}}(Q) \cap \mathbb{H}_Q.$  Since this geometry is always a half-ball, this constant is uniform.  Therefore, we have that for all $x \in B_{\frac{1}{4}}(Q) \cap \mathbb{H}_Q$
 \begin{align*}
 \phi(x) \le C L(x) \le C 2 C(n, \Lambda) \text{dist}(x, L_Q).
 \end{align*}
Thus, for $p \in \overline{T_{Q_0, r}\Omega} \cap B_{\frac{1}{4}}(Q)$, we have  
 \begin{align*}
T_{Q_0, r}u(p) \le \phi(p) \le C(n, \Lambda) \text{dist}(x, L_Q).
 \end{align*}
 Applying this argument to $\pm T_{Q_0, r}u$, we obtain the desired estimate.
\end{proof}

\begin{lem}\label{compactness 1}(Preliminary Compactness)
Let $u_i \in \mathcal{A}(n, \Lambda),$ $Q_i \in \partial \Omega_i \cap B_{1/4}(0)$, and $0<r_i \le \frac{1}{4}.$  Then there exists a subsequence  (also indexed by $i$) such that
\begin{enumerate}
\item $T_{Q_i, r_i}u_i \rightarrow u_{\infty}$ in $C^{0, \gamma}(\overline{B_1(0)}).$
\item If we define $\Omega_{\infty} = \emph{Interior}(\overline{\{|u_{\infty}|>0 \}})$, then $\overline{T_{Q_i, r_i}\Omega_i \cap B_{1}(0)} \rightarrow \overline{\Omega_{\infty} \cap B_{1}(0)}$ in the Hausdorff metric on compact subsets and $\Omega_{\infty}$ is a non-degenerate convex domain with $0 \in \partial \Omega_{\infty}$ which satisfies the same non-degeneracy as Corollary \ref{C: geometric non-degeneracy}.
\item $u_{\infty}$ is harmonic in $\Omega_{\infty}.$
\end{enumerate}
\end{lem}

\begin{proof}
By definition, $T_{Q_i, r_i}u_i(0) = 0$.  Therefore, Lemma \ref{L: unif Holder bound 1} implies the first convergence result by Arzela-Ascoli. Note that since $H_{T_{Q_i, r_i}\Omega_{i}}(0, 1, T_{Q_i, r_i}u_i) = 1$ for all $i$, $H_{\Omega_{\infty}}(0, 1, u_{\infty})=1$. 

By taking a further subsequence, we may assume that $\lim_i T_{Q_i, r_i}\Omega_i = \Omega'$ exists in a set theoretic sense.  The uniform convergence in (1) implies that for all $0< \epsilon$ $\{|u_{\infty}|> \epsilon\} \subset \bigcap_{N=1}^{\infty} \bigcup_{i=N}^{\infty}\{|T_{Q_i, r_i}u_i| \ge \epsilon - 1/N\}$.  Since $\{ |T_{Q_i, r_i}u_i| > 1/2\}$ is non-empty and $|T_{Q_i, r_i}u_i|> 1/4$ in $B_{(2C(n, \Lambda))^{-1/\gamma}}(\{ |T_{Q_i, r_i}u_i| > 1/2\})$, $\Omega_{\infty}$ is non-degenerate. For all $y \in \Omega_{\infty}$ and all $0< \delta$ there exists a point $x \in \Omega_{\infty}$ such that $|x - y| \le \delta$ and $|u_{\infty}(x)| > \epsilon$ for some $\epsilon>0$.  Since $x \in \overline{\lim}_iT_{Q_i, r_i}\Omega_i$ and $0<\delta$ was arbitrary, $\Omega_{\infty} \subset \lim_i B_{\delta}(T_{Q_i, r_i}\Omega_i)$ for all $\delta>0$. 

To see the converse containment, we observe that for any $x_0 \in \text{Int}(\Omega')$, if $B_{\delta}(x_0) \subset \Omega'$ then for all sufficiently large $i$ $B_{\delta/2}(x) \subset T_{Q_i, r_i}\Omega_i$ and $T_{Q_i, r_i}u_i$ is uniformly bounded in $W^{1, 2}(B_{\delta}(x_0))$. Therefore, $T_{Q_i, r_i}u_i$ converge to a harmonic function in $\text{Int}(\Omega')$. By the convexity of $T_{Q_i, r_i}\Omega_i$ for all $i$, $\text{Int}(\Omega')$ is connected. 

Now, let $0< \delta$ and suppose that for all sufficiently large $i$ there exists $x_i \in T_{Q_i, r_i}\Omega_i \setminus B_{\delta}(\Omega_{\infty}) \cap B_1(0)$. By passing to a subsequence, we may assume $x_i \rightarrow x_{\infty} \in B_1(0) \cap \overline{\Omega'} \setminus B_{\delta/2}(\Omega_{\infty})$
Since $u_{\infty}= 0$ on $B_{\delta/2}(x_{\infty})$ it must be that $\lim_{i \rightarrow \infty}\sup_{B_{\delta/2}(x_{\infty})}|T_{Q_i, r_i}u_i| = 0$. Since $text{Int}(\Omega')$ is connected, by unique continuation $T_{Q_i, r_i}u_i \rightarrow 0$ in $\Omega'$.  But this contradicts $u_{\infty}$ being non-trivial. Thus, there must be a subsequence such that $\overline{T_{Q_i, r_i}\Omega_i \cap B_{1}(0)} \rightarrow \overline{\Omega_{\infty} \cap B_{1}(0)}$ in the Hausdorff metric on compact subsets. Convexity and the conclusion of Corollary \ref{C: geometric non-degeneracy} are preserved under this mode of convergence, and so (2) is proved.  

Now that we know that $\Omega' = \Omega_{\infty}$, the previous argument proves (3), as well.
\end{proof}

\section{Estimates Off the Zero Set}\label{S: off zero set}

In this section we prove that analogs of the results of Sections \ref{sec:uniform frequency bound Singular set} and \ref{sec: compactness Zero Set} hold for all $p \in B_{\frac{1}{8}}(0) \cap \overline{\Omega}$ and $0< r \le \frac{1}{8}$. The key to obtaining estimates off $\{u=0\}$ is the following technical lemma.

\begin{lem}(Technical Lemma)\label{H lower bound}
For any $u \in \mathcal{A}(n, \Lambda)$, $Q \in B_\frac{1}{4}(0)$, and $0< r< \frac{1}{4},$ there is a constant $0< c(n, \Lambda)$ such that for all $y \in B_{1/2}(0),$ 
\begin{align*}
c(n , \Lambda) < H_{T_{Q, 2r}\Omega}(y, \frac{1}{4}, T_{Q, 2r}u) < C(n, \Lambda).
\end{align*}
\end{lem}

\begin{proof}
Note that the upper bound follows directly from Lemma \ref{L: unif Holder bound 1}.  To show the lower bound, we argue by compactness.  Suppose that there is a sequence of functions, $u_i \in \mathcal{A}(n, \Lambda)$, points $Q_i \in \partial \Omega_i \cap B_{\frac{1}{4}}(0)$ and radii $0< r_i< \frac{1}{4}$ such that there exist points $y_i \in B_{1/2}(0) \cap \overline{\Omega}$ for which 
\begin{align*}
H_{T_{Q_i, 2r_i}\Omega_i}(y_i, \frac{1}{4}, T_{Q_i, 2r_i}u_i) \le 2^{-i}
\end{align*}

Letting $i \rightarrow \infty$, by Lemma \ref{compactness 1}, there exists a subsequence $T_{Q_j, 2r_j}u_j$ which converges to a H\"older continuous function, $u_{\infty}$, which is harmonic in a non-degenerate convex domain, $\Omega_{\infty}$.  Note that $u_{\infty}$ vanishes on the boundary of $\partial \Omega_{\infty} \cap B_{8}(0)$.  Similarly, we may take subsequences such that $y_i \rightarrow y_{\infty}$.  Note that H\"older convergence implies $H_{\Omega_{\infty}}(0, 1, u_{\infty}) = 1$.
Since we have that $H_{\Omega_{\infty}}(y, \frac{1}{4}, u_{\infty}) = 0$, it must be that $u_{\infty} = u_{\infty}(y_{\infty})$ on $\partial B_{1/4}(y) \cap \Omega_{\infty}.$ If $\partial B_{1/4}(y) \subset \Omega_{\infty},$ then $u_{\infty} \equiv u_{\infty}(y_{\infty})$ in $\Omega_{\infty}.$  This contradicts $u_{\infty}(0) = 0$ and  $H_{\Omega_{\infty}}(0, 1, u_{\infty}) = 1$.  If $\partial B_{1/4}(y)$ intersects $\partial \Omega_{\infty}$, then $u_{\infty}(y_{\infty}) = 0$, since $u_{\infty}$ must vanish continuously on $\partial \Omega_{\infty}.$ However, this forces $u_{\infty} \equiv 0$, which contradicts  $H_{\Omega_{\infty}}(0, 1, u_{\infty}) = 1$. 
\end{proof}

\begin{lem}\label{N universal bound}(Bounding the Almgren Frequency)
Let $u \in \mathcal{A}(n, \Lambda)$, $p \in B_{\frac{1}{8}}(0) \cap \overline{\Omega}$ and $0< r \le \frac{1}{8}.$ Then, there is a constant, $C_2 = C_2(n, \Lambda)< \infty$ such that 
\begin{align*}
    N_{\Omega}(p, r, u) \le C_2.
\end{align*}
\end{lem}

\begin{proof}
For $p \in \{u=0\} \cap \overline{\Omega}$ Lemma \ref{N bound lem 1} proves the desired inequality. Let $p \in \Omega$ be such that $u(p) \not = 0$. Let $0< \delta = \text{dist}(p, \partial \Omega)$ and $Q \in \partial \Omega$ such that $|p-Q| = \delta.$  If $B_r(p) \subset \Omega$ we use the monotonicity of the Almgren frequency function to reduce to considering $B_\delta (p).$ We let $\tau = \max\{r, \delta\} \le \frac{1}{8}$

Consider $T_{Q, 4\tau}u$.  We note that
\begin{align*}
    C_1(n, \Lambda) & \ge N_{T_{Q, 4\tau}\Omega}(0, 1, T_{Q, 4\tau}u)\\
    & = \frac{\int_{B_1(0)}|\nabla T_{Q, 4\tau}u|^2 dV}{\int_{\partial B_1(0)}T_{Q, 4\tau}u^2d\sigma} \ge \int_{B_{\frac{1}{4}}(T_{Q, 4\tau}p)}|\nabla T_{Q, 4\tau}u|^2dV.
\end{align*}
On the other hand, Lemma \ref{H lower bound} implies that
\begin{align*}
    H_{T_{Q, 4\tau}\Omega}(T_{Q, 4\tau}p, \frac{1}{4}, T_{Q, 4\tau}u) & \ge c(n, \Lambda).
\end{align*}
Therefore, we have $N_{\Omega}(p, r, u) = \frac{D_{T_{Q, 4\tau}\Omega}(T_{Q, 4\tau}p, \frac{1}{4}, T_{Q, 4\tau}u)}{4 H_{T_{Q, 4\tau}\Omega}(T_{Q, 4\tau}p, \frac{1}{2}, T_{Q, 4\tau}u)} \le C(n, \Lambda).$
\end{proof}

\begin{lem}\label{Lipschitz bound}(Lipschitz Bounds)
For $u \in \mathcal{A}(n, \Lambda),$ for all $Q \in \partial \Omega \cap B_{\frac{1}{4}}(0)$ and all radii, $0< r \le \frac{1}{8}$, $T_{Q, r}u \in Lip(B_1(0))$ with uniform Lipschitz constant $Lip(T_{Q, r}u) \le C(n, \Lambda)$.
\end{lem}

\begin{proof}
Since $T_{Q, r}u$ is continuous and constant outside of $T_{Q, r} \Omega$, we reduce to bounding $\nabla T_{Q, r}u$ at interior points $y \in T_{Q, r}\Omega \cap B_1(0).$  Note that by our definition of the rescalings (Definition \ref{L2 rescaling def}) and Lemma \ref{L: H doubling-ish 1}
\begin{align*}
|T_{Q, r}u(y)| = & \left(\frac{\frac{1}{(4r)^{n-1}}H_{\Omega}(Q, 4r, u)}{\frac{1}{r^{n-1}}H_{\Omega}(Q, r, u)}\right)^{\frac{1}{2}} |T_{Q, 4r}u(y')|\\
\le & 4^{2C_1(n, \Lambda)} |T_{Q, 4r}u(y')|\\
\le & C(n, \Lambda) |T_{Q, 4r}u(y')|,
\end{align*}
where $y' = \frac{1}{4}y$.

Note that $y' \in B_{\frac{1}{4}}(0) \cap T_{Q, 4r}\Omega$.  Let $\delta = \text{dist}(y', T_{Q, 4r}\partial \Omega).$  Therefore, $\nabla T_{Q, 4r}u(y') = \fint_{B_{\delta}(y')} \nabla T_{Q, 4r}u dV$.  Recall that $|\nabla u|$ is subharmonic, and therefore by Lemma \ref{N universal bound}
\begin{align*}
|\nabla T_{Q, 4r}u(y')|  & \le \fint_{B_{\delta}(y')}|\nabla T_{Q, 4r} u|dV \le \left(\fint_{B_{\delta}(y')}|\nabla T_{Q, 4r} u|^2dV\right)^{\frac{1}{2}}\\
 & \le \left(C_2(n, \Lambda)\delta^{-2}\fint_{\partial B_{\delta}(y')}(T_{Q, 4r}u - T_{Q, 4r}u(y'))^2 d\sigma \right)^{\frac{1}{2}}.
\end{align*}

Now, let $Q' \in T_{Q, 4r} \partial \Omega$ be a point such that $\delta = |y' - Q'|$ and let $y' = y'' + Q'.$  Now, we translate the domain by $Q.$
\begin{align*}
\fint_{\partial B_{\delta}(y')}(T_{Q, 4r}u - T_{Q, 4r}u(y'))^2d\sigma =  \fint_{\partial B_{\delta}(y'')}(T_{Q, 4r}u(x + Q) - T_{Q, 4r}u(y''+ Q))^2d\sigma.
\end{align*}

Note that $T_{Q, 4r}u(x + Q') \in \mathcal{A}(n, C_1(n, \Lambda)).$  Now, by Corollary \ref{sub-linear bound} applied to $T_{Q, 4r}u(x + Q) \in \mathcal{A}(n, C_1(n, \Lambda))$  with $Q_0 = 0$ we bound 
\begin{align*}
\fint_{\partial B_{\delta}(y')}(T_{Q, 4r}u - T_{Q, 4r}u(y'))^2d\sigma & =  \fint_{\partial B_{\delta}(y'')}(T_{Q, 4r}u(x + Q) - T_{Q, 4r}u(y''+ Q))^2d\sigma(x)\\
& \le \fint_{\partial B_{\delta}(y'')}(4C(n, C_1(n, \Lambda))\delta)^2d\sigma  =  (4C(n, C_1(n, \Lambda))\delta)^2.
\end{align*}
Thus, we have that, 
\begin{align*}
|\nabla T_{Q, r}u(y)| \le & C(n, \Lambda) |T_{Q, 4r}u(y')|\\
& \le C(n, \Lambda)\left(C_1(n, \Lambda)\delta^{-2}\left(\fint_{\partial B_{\delta}(y)}(T_{Q, r}u - T_{Q, r}u(y))^2d\sigma \right)\right)^{\frac{1}{2}}\\
& \le C(n, \Lambda) C_1(n, \Lambda)^{\frac{1}{2}}\frac{1}{\delta}(4C(n, C_1(n, \Lambda))\delta) \le C(n, \Lambda).
\end{align*}
\end{proof}

We now prove that the Almgren frequency is a function of uniformly bounded variation.

\begin{lem}\label{bounded variation of N}(Bounded Variation)
Let $u \in \mathcal{A}(n, \Lambda)$ and $p \in B_{\frac{1}{8}}(0) \cap \overline{\Omega}$. Then, there is a constant $C_3 = C_3(n, \Lambda)< \infty$ such that for all $0< r \le \frac{1}{8}$
\begin{align*}
    \emph{var}(N_{\Omega}(p, r, u), [0, \frac{1}{8}]) \le C_3.
\end{align*}
\end{lem}

\begin{proof}
We estimate the variation by a ``rays of the sun" argument. Since $N_{\Omega}(p, r, u)$ is monotone increasing and bounded for $p \in \{u=0\} \cap \overline{\Omega}$, we argue for $u(p) \not = 0.$  Again, we let $\delta = \text{dist}(p, \partial \Omega).$
\begin{align*}
    \text{var}(N_{\Omega}(p, r, u), [0, \frac{1}{8}]) & \le 2\int_0^{\frac{1}{8}} N'_3(r)dr + |N_{\Omega}(p, 0^{+}, u) - N_{\Omega}(p, \frac{1}{8}, u)|\\
    & \le 2\int_{\delta}^{\frac{1}{8}} 2N_{\Omega}(p, r, u) \frac{1}{H_{\Omega}(p, r, u)}\int_{B_r(p) \cap \partial \Omega}u(p)\nabla u \cdot \vec \eta d\sigma dr + 2C_2.
\end{align*}
Now, if we let $Q_0 \in \partial \Omega$ be a point such that $|p - Q_0| = \delta$ we may calculate by Lemma \ref{H lower bound}, Lemma \ref{L: H doubling-ish 1}, Lemma \ref{sub-linear bound}, and Lemma \ref{Lipschitz bound}
\begin{align*}
\frac{1}{H_{\Omega}(p, r, u)}\int_{B_r(p) \cap \partial \Omega}u(p)\nabla u \cdot \vec \eta d\sigma & \le C(n, \Lambda)\int_{B_1(T_{Q_0, r}p) \cap T_{Q_0, r}\partial \Omega}T_{Q_0, r}u(y)\frac{1}{r}\nabla T_{Q_0, r}u \cdot \vec \eta_{T_{Q_0, r}\Omega} d\sigma\\
& \le C(n, \Lambda)\int_{B_1(T_{Q_0, r}p) \cap T_{Q_0, r}\partial \Omega}\frac{\delta}{r^2}d\sigma \le C(n, \Lambda)\frac{\delta}{r^2}.
\end{align*}
Thus, by Lemma \ref{N universal bound}, we may bound
\begin{align*}
    \text{var}(N_{\Omega}(p, r, u), [0, \frac{1}{8}]) & \le 8 C_2(n, \Lambda)C(n, \Lambda)\int_{\delta}^{\frac{1}{8}} \frac{\delta}{r^2}dr + 2C_2(n, \Lambda)\\
    & \le C(n, \Lambda) \delta(8 + \frac{1}{\delta}) + 2C_2(n, \Lambda) \le C_3(n, \Lambda).
\end{align*}
This proves the lemma.
\end{proof}

\begin{lem}(Compactness)\label{compactness}\label{strong convergence}
Let $u_i \in \mathcal{A}(n, \Lambda)$, $p_i \in \overline \Omega_i \cap B_{\frac{1}{8}}(0)$, and $r_i \in (0,\frac{1}{8}]$.  Then, there exists a subsequence and a function, $u_{\infty} \in W^{1, 2}_{loc}(\RR^n)$, such that $T_{p_j, r_j}u_j$ converges to $u_{\infty}$ in the following senses.
\begin{enumerate}
\item $T_{p_i, r_i}u_i \rightarrow u_{\infty}$ in $C^{0}(B_1(0))$.
\item $T_{p_i, r_i}u \rightarrow u_{\infty}$ in $L^2(B_1(0))$.
\item $\overline{T_{p_i, r_i}\Omega_i} \cap B_1(0) \rightarrow \overline{\Omega_{\infty}} \cap B_1(0)$ in the Hausdorff metric on compact subset and  $\overline{\Omega_{\infty}} \cap B_1(0) = \text{supp}(u_{\infty}) \cap B_1(0)$ is a non-degenerate, convex set.
\item $\nabla T_{p_i, r_i}u \rightarrow \nabla u_{\infty}$ in $L^2(B_1(0); \RR^n)$.
\end{enumerate}
\end{lem}

\begin{proof}  To see $(1)$, we observe that $T_{p_i, r_i}u_i(0) = 0$ and $\{T_{p_i, r_i}u_i\}$ are uniformly Lipschitz.  Therefore, by Arzela-Ascoli, there exists a subsequence which converges in $C^{0}(B_1(0))$. Since $C^{0}(B_1(0)) \subset L^2(B_1(0))$, this also proves $(2)$.

(3) follows analogously as in the proof of Lemma \ref{compactness 1}(2).  That is, if $\overline{T_{p_i, r_i}\Omega_i} \cap B_1(0) \not =  B_1(0)$ for a subsequence of $i$, then $\overline{T_{p_i, r_i}\Omega_i}$ is a translation of $\overline{T_{Q_i, r_i}\Omega_i}$ for some $Q_i \in \partial T_{p_i, r_i}\Omega_i \cap B_1(0)$. Thus, after possibly passing to a subsequence so that the $\lim_iQ_i$ exists, the argument of Lemma \ref{compactness 1}(2) applies.  (3) follows immediately.

By our choice of rescaling, $T_{p_j, r_j}u$, we have that $N_{\Omega}(0, 1, T_{p_j, r_j}u_j) = \int_{B_1(0)}|\nabla T_{p_j, r_j}u_j|^2dV.$  Therefore, Lemma \ref{N universal bound} gives that $\nabla T_{p_j, r_j}u_j$ are uniformly bounded in $L^2(B_1(0); \RR^n).$  Therefore, Rellich compactness gives weak convergence.

The only thing remaining to show is that $\nabla T_{p_j, r_j}u_j \rightarrow \nabla u_{\infty}$. By (3), we may choose a subsequence such that, $\partial \Omega_j$ have a convergent subsequence such that $T_{p_i, r_i}\partial \Omega_i \rightarrow \partial \Omega_{\infty}$ locally in the Hausdorff metric to a non-degenerate convex domain.  Since the boundary of a convex domain is locally the graph of a Lipschitz function $\overline{\dim}_\mathcal{M}(\partial \Omega_\infty \cap B_1(0)) = n-1$.  Thus, by continuity of measures and Lemma \ref{Lipschitz bound}, for all $\epsilon > 0$ we can find a $\tau(\Lambda, n, \epsilon)$ independent of $T_{p_j, r_j}u_j$, such that
\begin{align*}
\int_{B_1(0) \cap B_{\tau}(\partial \Omega_{\infty})} |\nabla T_{p_j, r_j}u_j|^2dV \le \epsilon.
\end{align*}
Therefore, using the notation $\partial \Omega_{j, \tau} = B_{\tau}(T_{p_j, r_j}\partial \Omega_{j})$
\begin{align*}
\lim_{j \rightarrow \infty} D_{\Omega_i}(1, 0, T_{p, r_j}u_j)  = & \lim_{j \rightarrow \infty} \int_{B_1(0)}|\nabla T_{p_j, r_j}u_j|^2dV\\
 = & \int_{B_1(0) \setminus \partial \Omega_{j, \tau}} |\nabla T_{p_j, r_j}u_j|^2dV + \lim_{j \rightarrow \infty} \int_{B_1(0) \cap \partial \Omega_{j, \tau}}|\nabla T_{p, r_j}u_j|^2dV \\
 \le & \lim_{j \rightarrow \infty} \int_{B_1(0) \setminus B_{\tau/2}(\partial \Omega_{\infty})} |\nabla T_{p_j, r_j}u_j|^2dV + \epsilon \\
 \le &  D_{\Omega_{\infty}}(1, 0, u_{\infty}) + \epsilon,
\end{align*}
where the last equality follows from $W^{1, 2}$-convergence of harmonic functions in the region $B_1(0) \setminus B_{\tau}( \partial \Omega_{\infty}).$  Since $\epsilon > 0$ was arbitrary, we have that $\lim_{j \rightarrow \infty}D_{\Omega_i}(1, 0, T_{p_j, r_j}u_j) \le D_{\Omega_{\infty}}(1, 0, u_{\infty})$.  The other inequality follows from the same trick or from lower semi-continuity. Thus, $\lim_{j \rightarrow \infty}D(1, 0, T_{Q_j, r_j}u_j) = D_{\Omega_{\infty}}(1, 0, u_{\infty})$.  This implies strong convergence.
\end{proof}

\begin{cor}(Convergence of the Almgren frequency)\label{N continuity}
For $u_j \in \mathcal{A}(n, \Lambda)$, $p_j \in B_{\frac{1}{8}}(0) \cap \overline \Omega_i$, and $r_j \in (0, \frac{1}{8}]$, there exists a subsequence and a limit function such that 
\begin{align}
    N_{T_{p_j, r_j}\Omega_i} (0, 1, T_{p_j, r_j}u_j) \rightarrow N_{\Omega_{\infty}}(0, 1, u_{\infty}).
\end{align}
\end{cor} 

\begin{proof}  
The continuous convergence of $T_{p_j, 2r_j}u_j$ in $B_{1}(0)$ and the strong convergence $\nabla T_{p_j, 2r_j}u_j$ in $B_{1}(0)$ give the desired convergence of $H_{T_{p_j, 2r_j}\Omega_i} (0, \frac{1}{2}, T_{p_j, 2r_j}u_j)$ and $D_{T_{p_j, 2r_j}\Omega_j} (0, \frac{1}{2}, T_{p_j, 2r_j}u_j)$, respectively.  Recall that by Definition \ref{N def} and Definition \ref{L2 rescaling def} 
\begin{align*}
N_{T_{p, 2r}\Omega}(0, \frac{1}{2}, T_{p, 2r}u) & = \frac{1}{2}\frac{D_{T_{p, 2r}\Omega}(0, \frac{1}{2}, T_{p, 2r}u)}{H_{T_{p, 2r}\Omega}(0, \frac{1}{2}, T_{p, 2r}u)}\\
& = \frac{D_{T_{Q, r}\Omega}(0, 1, T_{Q, r}u)}{H_{T_{Q, r}\Omega}(0, 1, T_{Q, r}u)}.
\end{align*}
\end{proof}

\begin{cor}(Limit functions are harmonic in the limit domain)\label{C: limit harmonic}
Let the sequence of functions $T_{p_j, r_j}u$ converge to the function $u_{\infty}$ in the senses of Lemma \ref{compactness}.  Then, $u_{\infty}$ is harmonic in $\Omega_{\infty}$.
\end{cor}

\begin{proof}  Recall that the boundaries, $T_{Q, r}\partial \Omega_j \rightarrow \partial \Omega_{\infty}$ in the Hausdorff distance on compact subsets.  Therefore, for any $0 < \epsilon$, for $j$ large enough, every $T_{Q, r}u$ will be harmonic in the region $B_1(0) \setminus B_{\epsilon}(\partial \Omega_{\infty}).$   By $C^{0, \gamma}(B_1(0))$ convergence of harmonic functions, $u_{\infty}$ is therefore harmonic in $B_1(0) \setminus B_{\epsilon}(\partial \Omega_{\infty}).$  Letting $\epsilon \rightarrow 0$ gives the desired statement.
\end{proof}

\section{Geometric Control}\label{S: quant rigidity}

The main results of this section are two ``quantitative rigidity" results about homogeneous harmonic functions.  Both are essentially consequences of the compactness obtained in Lemma \ref{strong convergence}. 

\begin{lem}\label{N constant gives homogeneity}
Let $u \in \mathcal{A}(n, \Lambda).$  Let $p \in \overline{\Omega} \cap B_{\frac{1}{8}}(0)$ and $0< r \le \frac{1}{8}.$ If
\begin{align*}
N_{\Omega}(p, r, u)= N_{\Omega}(p, r/10, u),
\end{align*}
and either $u(p)=0$ or $B_r(p) \subset \Omega$ then $u$ is $(0, 0, s, p)$-symmetric for all $0<s\le 1$. 
\end{lem}

\begin{proof}
The hypotheses imply that $N_{\Omega}(p, s, u)$ is a constant for all $r/10 \le s \le r$.  Furthermore, using the notation in Lemma \ref{N derivative calculation},
 $N'_1(s) = N'_2(s) = N'_3(s) = N'_4(s)=0$ for all $r/10\le s \le r.$  Thus, by Lemma \ref{N_1 redefine} we have that for all $y \in \partial B_s(p)$
\begin{align*}
 \nabla u \cdot (y-p) = & \partial_{r} u \\
= &  N_{\Omega}(p, s, u) (u - u(p)).
 \end{align*}
Since $N_{\Omega}(p, s, u)$ is a constant for all $r/10\le s\le r$, this becomes a separable ODE in polar coordinates, and $u -u(p) = v(s, \theta) = s^{N_{\Omega}(p, r, u)}v(\theta).$  Since $\Omega \cap B_r(p)$ is open, unique continuation implies that $u-u(p)$ is a homogeneous function of degree $N_{\Omega}(p, r, u)$ in $\Omega \cap B_1(p).$ In particular, $u$ is $(0,0,s, p)$-symmetric for all $0<s \le 1$.  

\end{proof}

Standard quantitative rigidity results usually prove that if $N_{\Omega}(p, r, u)$ is \textit{almost} constant ($N_{\Omega}(p, 1, u) - N_{\Omega}(p,  r, u) \le \delta(\epsilon) $), then $u$ is \textit{almost} a homogeneous harmonic polynomial ($||T_{p, 1}u - P||_{L^2(B_1(0))} \le \sqrt{\epsilon}$).  However, these results rely essentially upon the monotonicity of $N_{\Omega}(p, r, u)$.  If $N_{\Omega}(p, r, u)$ is not monotonic, then $N_{\Omega}(p, 1, u) = N_{\Omega}(p, r, u)$ does not imply that the Almgren frequency is constant.  In fact, even if $N_{\Omega}(p, r, u)$ is constant for $1/10 < r< 1$, if $u(p) \not = 0$, it is not clear that $u$ would be homogeneous. To overcome this technical issue we consider $p$ which are merely very close to $\{u=0\}$. 

\begin{lem} (Quantitative Rigidity)\label{quant rigidity}
Let $u \in \mathcal{A}(n, \Lambda)$, as above.  Let $p \in B_{\frac{1}{8}}(0) \cap \overline \Omega$ and $0< r \le \frac{1}{8} $.  For every $\delta > 0$, there is an $0< \gamma_0= \gamma_0(n, \Lambda, \delta)$ such that for any $0 < \gamma \le \gamma_0$ if
\begin{align*}
|N_{\Omega}(0, 1, T_{p, r}u) - N_{\Omega}(0, 1/10, T_{p, r}u)| \le \gamma
\end{align*}
and either $\text{dist}(p, \{u=0\}) \le \gamma r$ or $B_1(0) \subset T_{p, r}\Omega$ then $T_{p, r}u$ is $(0, \delta, 1, 0)$-symmetric.
\end{lem}

\begin{proof}  We argue by contradiction.  Assume that there exists a $\delta >0$ such that there is a sequence of functions, $u_i \in \mathcal{A}(n, \Lambda)$, and points, $p_i \in B_{\frac{1}{8}}(0) \cap \overline \Omega_i$, radii $0< r_i < \frac{1}{8}$, such that $\text{dist}(p_i, \{u_i=0\}) \le r_i 2^{-i}$ and
$$
|N_{\Omega_i}(0, 1, T_{p_i, r_i}u_i) - N_{\Omega_i}(0, 1/10, T_{p_i, r_i}u_i)| \le 2^{-i}
$$
but that no $T_{p_i, r_i}u_i$ is $(0, \delta, 1, 0)$-symmetric.

By Lemma \ref{strong convergence} we have that there exists a subsequence such that $T_{p_j, r_j}u_j$ converges strongly in $W^{1, 2}(B_1(0))$ to a function $u_{\infty}$.  By Corollary \ref{C: limit harmonic}, we know that $u_{\infty}$ is harmonic in a convex domain $\Omega_{\infty}$. Furthermore, by Lemma \ref{strong convergence}(1) $u_{\infty}(0) = 0$ and $u_{\infty} = 0$ on $\partial \Omega_{\infty}$. By Lemma \ref{N universal bound} and the proof of Corollary \ref{N continuity} applied to $T_{p_j, r_j}u_j$, we have that $\lim_{j \rightarrow \infty} N_{\Omega_j}(0, r, T_{p_j, r_j}u_j) = N_{\Omega_{\infty}}(0, r, u_{\infty}) \in [0, C_2(n, \Lambda)]$, and $N_{\Omega_{\infty}}(0, r, u_{\infty})$ is constant for $1/10 \le r \le 1$.  Therefore, by Lemma \ref{N constant gives homogeneity} $u_{\infty}$ is $(0, 0, 1, 0)$-symmetric. This contradicts our assumption that that no $T_{p_i, r_i}u_i$ is $(0, \delta, 1, 0)$-symmetric.

The case that $B_{1}(0)\subset T_{p_i, r_i}\Omega_i$, repeat the argument to obtain the same contradiction. 
\end{proof}

\begin{rmk}\label{quant rigidity invariances}
By the scale invariance of the Almgren frequency, Lemma \ref{quant rigidity} implies that for all $0< \delta$, if $0< \gamma \le \gamma_0(n, \Lambda, \delta)$, then 
\begin{align*}
|N_{\Omega}(p, r, u) - N_{\Omega}(p, r/10, u)| \le \gamma
\end{align*}
and either $\text{dist}(p, \{u=0\}) \le \gamma r$ or $B_r(p) \subset \Omega$ implies $u$ is $(0, \delta, r, p)$-symmetric.
\end{rmk}

Next, we obtain a ``cone-splitting" result. The prototypical example of a result like this is the following proposition.  See \cite{HanLin_nodalsets} Theorem 4.1.3 for the proof of similar results.

\begin{prop}\label{P: classic cone splitting}
Let $P: \RR^n \rightarrow \RR$ be a $0$-symmetric function.  Let $k \le n-2$.  If $P$ is symmetric with respect to some $k$-dimensional subspace $V$ and $P$ is homogeneous with respect to some point $x \not\in V$, then $P$ is $(k + 1)$-symmetric with respect to $span\{x, V\}$.  
\end{prop}

We now prove a similar result for our almost-symmetric functions $u \in \mathcal{A}(n, \Lambda)$.

\begin{lem}\label{cone splitting}(Cone-splitting)
Let $0\le k \le n-2$ be an integer and $u \in \mathcal{A}(n, \Lambda)$.  For any fixed $\epsilon, \tau > 0$ and any $0< r \le 1/8$, there is a $0< \delta_0 = \delta_0(n, \tau, \epsilon, \Lambda)$ such that for $0< \delta< \delta_0$, the following holds.  If  $p \in \overline{\Omega} \cap B_{\frac{1}{8}}(0)$ and $u$ is $(k, \delta, r, p)$-symmetric with respect to a $k$-dimensional subspace $V$ and $(0, \delta, r, x)$-symmetric for some $x \in B_r(p)\setminus B_{\tau r}(V+p)$, then $v$ is $(k+1, \epsilon, r, p)$-symmetric.
\end{lem}

\begin{proof} Assume that there exists a $\delta, \tau > 0$ for which there exist a sequence of $0< r_i \le 1$, function, $u_i \in \mathcal{A}(n, \Lambda)$ and points $\{p_i\}$ for which $u_i$ is $(k, i^{-1}, r_i, p_i)$-symmetric with respect to some $V_i$ and $(0, i^{-1}, r_i, x_i)$-symmetric for some $x_i \in B_{r_i}(p_i) \setminus B_{\tau}(V_i+p_i)$, but that all $u_i$ are not $(k+1, \delta, r_i, p_i)$-symmetric.

By considering $T_{p_i, r_i}u_i$ and applying Lemma \ref{strong convergence} and Lemma \ref{C: limit harmonic} there exists a function $u_{\infty} \in L^2(B_1(0))$ such that a subsequence $T_{p_j, r_j}u_j \rightarrow u_{\infty}$ in the senses of the lemma. Note that $u_{\infty}$ is non-degenerate. Taking further subsequences, we may reduce to a sequence for which 
$$
T_{p_j, r_i}V_j \rightarrow V, \quad x_j \rightarrow x \in \overline{B_1(0)}\setminus B_{\tau}(V).$$  
Note that $u_{\infty}$ is  $(k, 0, 1, 0)$-symmetric with respect to $V$, and harmonic. Therefore, $u_{\infty}$ is $(k+1, 0, 1, 0)$-symmetric.  Since $u_i \rightarrow u$ in $L^2(B_2(0))$, we have our contradiction. By taking the smallest $0< \delta_0$ for $0 \le k \le n-2$, we eliminate the dependence upon $k$.   
\end{proof}

\begin{cor}\label{c: application of cone splitting}
Let $k \le n-2$ and $u \in \mathcal{A}(n, \Lambda)$.  For any fixed $\epsilon, \tau > 0$ and any $0< r \le 1$, there is a $0< \delta_0 = \delta_0(n, \tau, \epsilon, \Lambda)$ such that for $0< \delta< \delta_0$, the following holds.  If  $p \in \overline{\Omega} \cap B_{\frac{1}{4}}(0)$ and $u$ is $(0, \delta, r, p)$-symmetric but not $(k+1, \epsilon, r, p)$-symmetric then there exists an affine $k$-plane $V$ such that
\begin{align}
    \{x \in B_r(p): u \text{ is } (0, \delta, r, x)\text{-symmetic} \} \subset B_{\tau r}(V).
\end{align}
\end{cor}
\begin{proof}
Let $0< \epsilon_1, ..., \epsilon_{n-1}$ be small parameters to be chosen later.  To prove the lemma, we inductively apply Lemma \ref{cone splitting}.  That is, suppose that there was another point $x_1 \in B_r(p) \setminus B_{\tau r}(p)$ such that $u$ was $(0, \delta, r, x)$-symmetic, then Lemma \ref{cone splitting} implies that for any $0< \epsilon_1$ there is a $0< \delta_0(n, \tau, \epsilon_1, \Lambda)$ such that if $\delta \le \delta_0(n, \tau, \epsilon_1, \Lambda)$ then $u$ is $(1, \epsilon_1, r, p)$-symmetric with respect to some affine $1$-plane $V_1$. Inducting up, assume that $u$ is $(i, \epsilon_{i}, r, p)$-symmetric with respect to $V_i$.  Then, for any $0<\epsilon_{i+1}$ there exists a $\delta_0(n, \tau, \epsilon_{i+1}, \Lambda)$ such that if we can find an $x_{i+1} \in B_r(p) \setminus B_{\tau r}(V_{i})$ such that $u$ is $(0, \delta, r, x)$-symmetric and $0< \delta, \epsilon_i \le \delta_0(n, \tau, \epsilon_{i+1}, \Lambda)$, then $u$ is $(i+1, \epsilon_{i+1}, r, p)$-symmetric with respect to $V_i$.

Therefore, choosing $0< \epsilon_{n-1} \le \delta_0(n, \tau, \epsilon, \Lambda)$ and $0< \epsilon_{i-1} \le \delta_0(n, \tau, \epsilon_i, \Lambda)$ for all $i = 1, ...., n-1$, we obtain $0< \delta_0(n, \tau, \epsilon, \Lambda) = \delta_0(n, \tau, \epsilon_1, \Lambda)$ as desired.  Since it is assumed that such that $u$ is not $(k+1, \epsilon, r, p)$-symmetric, this procedure must terminate before step $k+1$ and therefore there must exist an affine $k$-plane such that the claim of the corollary holds.
\end{proof}


\section{The Covering and its Properties}\label{S:covering and properties}
The lemmata in the previous section allow us to inductively define a covering with the right packing conditions.  Quantitative rigidity allows us to prove a ``Quantitative Differentiation" lemma that bounds the number of scales across which the frequency can change by more than some threshold $\gamma>0$.  Cone splitting on the other hand, will give us good geometric control of the singular set at scales for which $v$ is close to a homogeneous harmonic polynomial.  Together, these things will give us the necessary packing conditions.  First, we describe the covering.

\subsection{The General Construction}\label{S:the construction}

Let $u \in \mathcal{A}(n, \Lambda)$, and let $\epsilon, r >0$, $k \le n-2$, and $N \in \mathbb{N}$ be given.  We use the notation $\rho_i = 10^{-i}$. In this section we describe a general procedure which will produce a cover of $\mathcal{C}^k_{\epsilon, r}(u) \cap B_{\frac{1}{10}}(0)$ by balls of radius $\rho_N$.

We begin by defining an auxiliary quantity.  Let 
\begin{align}
    \mathcal{D}(u, x, r) = \inf\{\delta' > 0 : u \text{  is  } (0, \delta', r, x)\text{-symmetric}\}.
\end{align}
Let $0< \delta_0$. We shall refer to $\delta_0$ as the sorting threshold.  For any $i \in \mathbb{N}$ we can assign to each $x \in \mathcal{C}^k_{\epsilon, r}(u) \cap B_{\frac{1}{8}}(0)$ an $i$-tuple $T^i(x)$ according to the rule 
\begin{align*}
    (T^i(x))_j & = 1 \qquad \text{if  } \mathcal{D}(u, x, \rho_j) \ge \delta_0\\ (T^i(x))_j & = 0 \qquad \text{if  } \mathcal{D}(u, x, \rho_j) < \delta_0.
\end{align*}
For any $T_i$ we shall use $|T^i|$ to denote the sum of the entries.  Note that there is a partial ordering on the set of these $i$-tuples.  That is, if $k < i$, we can say that $T^k < T^i$ if $(T^k)_j= (T^i)_j$ for all $j \in \{1, 2, ..., k \}$.

Now, we partition our set according to these $i$-tuples.  For any given $i$-tuple, $T^i \in \{0, 1 \}^i$, we define
\begin{equation*}
E(T^i) = \{x \in \mathcal{C}^k_{\epsilon, r}(u) \cap B_{\frac{1}{10}}(0) : T^i(x) = T^i \}.
\end{equation*}
It follows immediately from the definitions that $E(T^i) \subset E(T^k)$ if and only if $T^k < T^i$.

We now define our covering inductively.  For $i=1$, we let $C^k_{\epsilon, r}(T^i) = B_{\frac{1}{10}}(0)$ for both $1$-tuples $T^i \in \{0, 1\}^1$.  Now, assume that $i \in \mathbb{N}$, $i < N$, and $C^k_{\epsilon, r}(T^i)$ has been defined and consists of balls of radius $\rho_{i}$.  Within each ball $B_{\rho_{i}}(y) \in C^k_{\epsilon, r}(T^i)$ partition the set $B_{\rho_{i}}(y) \cap E(T^{i})$ into the sets $E(T^{i+1})$ for $T^{i+1}$ such that $T^{i} < T^{i+1}$. For either such $T^{i+1}$, take a minimal covering of $B_{\rho_{i}}(y) \cap E(T^{i+1})$ by balls of radius $\rho^{i+1}$ centered at points in $B_{\rho_{i}}(y) \cap E(T^{i+1})$.  The union of these balls is $C^k_{\epsilon, r}(T^{i+1})$.

For some $i$-tuples, the set $E(T^i)$ may be empty.  In this case, we simply allow the corresponding collection of balls, $C^k_{\epsilon, r}(T^i)$, be empty. 

If $i = N$, we terminate the procedure. Note that for any sorting threshold $0< \delta_0$ and $N \in \mathbb{N}$ this procedure defines a sequence of collections such that
\begin{align*}
    \mathcal{C}^k_{\epsilon, r}(u) \cap B_{\frac{1}{10}}(0) \subset \bigcup_{T^N} \bigcup_{B_{\rho_N}(y) \in C^k_{\epsilon, r}(T^N)} B_{\rho_N}(y).
\end{align*}

\subsection{Properties of the Construction}

Now, we argue that there is a choice of sorting threshold $0<\delta_0$ with the desired properties.

\begin{lem}\label{bounded number of collections in the cover}
Let $u \in \mathcal{A}(n, \Lambda)$.  Let $0< \delta'$ be the sorting threshold in the construction above.  For any $i \in \mathbb{N}$ there are at most $N^{D(n,  \Lambda, \delta')}$ nonempty sets $E(T^i)$ such that $E(T^i)$ is non-empty.
\end{lem}

\begin{proof} Let $0 < \delta'$ be given.  Let $0<\gamma_0(n, \Lambda, \delta')<<1$ as in Lemma \ref{quant rigidity} and Remark \ref{quant rigidity invariances}. Now, decompose $\Omega = \cup_{j=0}^{\infty}A_j(\Omega)$ where 
\begin{align*}
    A_j(\Omega) = \Omega \cap B_{\gamma_0^{j}}(\partial \Omega) \setminus B_{\gamma_0^{j+1}}(\partial \Omega).
\end{align*}
We shall argue that there is a $D = D(n, \Lambda, \delta')$ such that $|T^N(p)| \le D$ for all $p \in \overline{\Omega} \cap B_{10}(0)$.  If the claim is true, then if $N \le D$ there are at most $2^N \le N^D$ $N$-tuples with $|T^N| \le D$. And, if $N \ge D$ there are at most $\binom{N}{D}$ many $N$-tuples with $|T^N| \le D$.  Since $\binom{N}{D}$ $\le N^D$ we have the desired claim: $\mathcal{C}^k_{\epsilon, r}(u) \cap B_{1/10}(0)$ is contained in the union of at most $N^D$ nonempty sets $E(T^N)$ and covered by at most $N^D$ collections $C^k_{\epsilon, r}(T^N)$.

To prove the claim, we argue by cases.  Let $p \in E(T^N)$. If $p \in \partial \Omega$, then $u(p)=0$ and 
\begin{align*}
    |T^N(p)| \le \left|\{i \in \mathbb{N} : |N_{\Omega}(p, \rho_i, u) - N_{\Omega}(p, \rho_{i+1}, u)| \ge \gamma_0(n, \Lambda, \delta') \}\right| \le C_3(n, \Lambda)\gamma_0^{-1}.
\end{align*}
for $\gamma_0(n, \Lambda, \delta')$ as in Remark \ref{quant rigidity invariances}.

Now, assume that $p \in A_j(\Omega)$. Note that by Remark \ref{quant rigidity invariances},  
\begin{align*}
    |N_{\Omega}(p, \rho_i, u) - N_{\Omega}(p, \rho_{i+1}, u)|\le \gamma_0(n, \Lambda, \delta')
\end{align*}
still implies $(T^N(p))_i = 0$ if $B_{\rho_i}(p) \subset \Omega$ or if $\text{dist}(p, \partial \Omega) \le \gamma_0(n, \Lambda, \delta')\rho_i$. 
But, the conditions $B_{\rho_i}(p) \subset \Omega$ or $\text{dist}(p, \partial \Omega) \le \gamma_0(n, \Lambda, \delta')\rho_i$ only fail for $i \in [(j-1)\frac{\ln(\gamma_0)}{\ln(\rho_0)}, (j+1)\frac{\ln(\gamma_0)}{\ln(\rho_0)}]$.  Therefore 
\begin{align*}
    |T^N(p)| \le \frac{C_3(n, \Lambda)}{\gamma_0(n, \Lambda, \delta')} + 3\frac{\ln(\gamma_0(n, \Lambda, \delta'))}{\ln(\rho_0)} = D(n, \Lambda, \delta')
\end{align*} independent of $j$. This proves the claim. 
\end{proof}

\begin{rmk}\label{r: new gamma}
If, later, we choose $0<\gamma \le \gamma_0(n, \Lambda, \delta')$ so that $\gamma= \gamma(n, \Lambda, \epsilon)$, then the statement of Lemma \ref{bounded number of collections in the cover} holds with a new $D = D(n, \Lambda, \epsilon)$.
\end{rmk}

We now prove that this construction satisfies the claimed packing condition.

\begin{lem}\label{construction packing properties}
Let $u \in \mathcal{A}(n, \Lambda)$, $0< \epsilon, r <1$, and $k \le n-2$. Let $\delta' = \delta_0(n, 1/10, \epsilon, \Lambda)$ as in Corollary \ref{c: application of cone splitting}.  Let $\gamma_0(n, \Lambda, \delta')$ be as in Lemma \ref{quant rigidity}. 

Then for all $N \in \mathbb{N}$ with $\rho_N \ge r$ there exist constants $0< c_1, c_2$ depending only on the ambient dimension $n$ and a constant $D(n, \Lambda, \epsilon)$ such that each collection $C^k_{\epsilon, r}(T^N)$ consists of at most $(c_1\rho_1^{-n})^D (c_2\rho_1^{-k})^{N-D}$ balls of radius $\rho_N$.
\end{lem}

\begin{proof} For any given $N$-tuple $T^N$ for which $E(T^N)$ is non-empty, let $T^{i} < T^N$. Now, let $B_{\rho_{i}}(x) \in C^k_{\epsilon, r}(T^{i})$. Consider the set $$A := B_{\rho_{i}}(x) \cap E(T^N).$$  We argue by cases. 

\textit{Case 1.} $(T^N)_{i} = 0$. In this case, $u$ is $(0, \delta_0, \rho_i, x)$-symmetric. By applying Corollary \ref{c: application of cone splitting}, we see that $A \subset B_{\rho_{i}}(x) \cap B_{\rho_{i+1}}(V^k)$ for some $k$-dimensional plane $V^k$. Thus, the minimal covering from the construction can cover $A \subset B_{\rho_{i}}(x) \cap B_{\rho_{i+1}}(V^k)$ by at most $c_1(n)\rho_1^{-k}$ balls of radius $\rho_{i+1}$.

\textit{Case 2.}  $(T^N)_{i} = 1$.  In this case, we have no control. Therefore, the minimal covering of $A$ described in the construction could consist of at most $c_2(n)\rho_1^{-n}$ balls of radius $\rho_{i+1}$.

Carrying this process through, by the proof of Lemma \ref{bounded number of collections in the cover} Case 2 can only happen at most $D(n, \Lambda, \epsilon)$ times.  Thus, $C^k_{\epsilon, r}(T^N)$ is a collection of at most $(c_1\rho_1^{-n})^{D} (c_2\rho_1^{-k})^{N-D}$ balls of radius $\rho^N$, as claimed.
\end{proof}

\subsection{Proof of Theorem \ref{T: main theorem 1}}

\begin{proof}  Let $0< \epsilon, r <1$, and $k \le n-2$. Let $\delta_0 = \delta_0(n, 1/10, \epsilon, \Lambda)$ as in Corollary \ref{c: application of cone splitting}.  Let $0<\gamma_0(n, \Lambda, \delta_0)$ be as in Lemma \ref{quant rigidity}. Recalling $c_2(n)$ from Lemma \ref{construction packing properties}, let
\begin{align*}
    0<\gamma \le \min\{\gamma_0(n, \Lambda, \delta_0), c_2^{\frac{-2}{\epsilon}}\} < 1.
\end{align*}

The construction given in Section \ref{S:the construction} gives a covering of $B_{\frac{1}{10}}(0) \cap \mathcal{C}^k_{\epsilon, r}(u)$ by balls of radius $\rho_N$.  Doubling the radius of these balls is sufficient, then, to cover $B_{\frac{1}{10}}(0) \cap B_{\rho^N}(\mathcal{C}^k_{\epsilon, r}(u))$.

Thus, by Remark \ref{r: new gamma} for $D = D(n, \Lambda, \epsilon)$ 
\begin{align*}
\text{Vol}(B_{\frac{1}{10}}(0) \cap B_{\rho_N}(\mathcal{C}^k_{\epsilon, r}(u))) & \le  N^D(c_1\rho_1^{-n})^D (c_2\rho_1^{-k})^{N-D}(\omega_n 2\rho^{N})^n\\
& \le  c(n) N^D c_1^D c_2^{N-D} (\omega_n2)^n \rho_N^{n-k} \rho_1^{-D(n-k)}.
\end{align*}

Estimating $\rho_1^{-D(n-k)} \le C(n, \Lambda, \epsilon)$ and $N^D \le C(n, \Lambda, \epsilon) c_2(n)^N$ for all $N \in \mathbb{N}$ and recalling from our choice of $0<\gamma$ that $c_2(n) \le \gamma^{-\frac{\epsilon}{2}}$, we obtain
\begin{eqnarray*}
\text{Vol}(B_{\frac{1}{8}}(0) \cap B_{\rho_N}(\mathcal{C}^k_{\epsilon, r}(u))) & \le & N^D c_1^D c_2^{N-D} (\omega_n2)^n \rho_N^{n-k} \rho_1^{-D(n-k)}\\
 & \le & C(n, \Lambda, \epsilon) c_2(n)^{2N-D}  c_1(n)^D (\omega_n2)^n \rho_N^{n-k}\\
& \le & C(n, \Lambda, \epsilon)  \gamma^{-N\epsilon}\rho_N^{n-k}\\
& \le & C(n, \Lambda, \epsilon) \rho_N^{n-k-\epsilon}.
\end{eqnarray*}
Thus, for any $0< r$ we may find an $N \in \mathbb{N}$ such that $\rho_{N+1} \le r < \rho_{N}$.  Thus, we may estimate
\begin{eqnarray*}
\text{Vol}(B_{\frac{1}{10}}(0) \cap B_{r}(\mathcal{C}^k_{\epsilon, r}(u))) 
& \le & C(n, \Lambda, \epsilon) \rho_{N}^{n-k-\epsilon}\\
& \le & C(n, \Lambda, \epsilon) \left(10 r\right)^{n-k-\epsilon} \le C'(n, \Lambda, \epsilon) \left(r\right)^{n-k-\epsilon}.
\end{eqnarray*}

Covering $B_\frac{1}{4}(0)$ $C(n)$ many such balls of radius $1/10$ and repeating the argument produces \eqref{e:main theorem 1 estimate}.
\end{proof}

\section{Containment}\label{S:containment}

We now turn to proving Lemma \ref{e-reg containment}, which follows from the rigidity and continuity of the Almgren frequency function.

\begin{lem}\label{force r small}
Let $u \in \mathcal{A}(n, \Lambda)$.  Suppose that $Q \in \partial\Omega \cap B_{1/4}(0)$ and there exists a $p \in \mathcal{C}^{n-2}(u)$ such that $B_r(p) \subset \Omega \cap B_{\frac{1}{4}}(0)$ and $p \subset B_{2r}(Q)$.  Then there is an $0<\delta(n, \Lambda)$ such that $N_{\Omega}(Q, 4r, u) \ge 1+ \delta$.
\end{lem}

\begin{proof}
Suppose that the statement is false.  Then, there exists a sequence of $u_i$, $Q_{i}, p_{i}$ and $0<r_i$ such that $B_{r_i}(p_{i}) \subset \Omega_{i} \cap B_\frac{1}{4}(0)$ for which \begin{align*}
    N_{T_{Q_{i},4r_i}\Omega_i}(0, 1, T_{Q_{i}, 4r_i}u_i) \le 1+ 2^{-i}.
\end{align*}
By Lemma \ref{compactness}, we may extract a subsequence $T_{Q_{j}, 4r_j}u_j$ which converges strongly in $W^{1, 2}(B_1(0))$ to a function $u_{\infty}.$ Similarly, we may assume that $T_{Q_{j}, 4r_j}p_j \rightarrow p$ and that $B_{1/4}(p) \subset \overline{\Omega_{\infty} \cap B_{\frac{1}{2}}(0)}.$  Observe that by Lemma \ref{N continuity}, $N_{\Omega_{\infty}}(0, 1, u_\infty) = 1$.  Therefore, $N_{\Omega_\infty}(0, r, u_\infty)$ is constant for $0< r \le 1$.  By Lemma \ref{N constant gives homogeneity}, $u_{\infty}$ must be a piecewise linear function.  Therefore, $N_{\Omega_\infty}(p, 1/4, u_{\infty}) = 1$. However, $N_{\Omega_j}(p_j, r_j, u_j) \ge 2$ for all $j$ which implies by Lemma \ref{N continuity} that $N_{\Omega_\infty}(p, 1/4,u_{\infty})\ge 2$. This is a contradiction.  
\end{proof}

\begin{lem}\label{boundary points not n-1 symmetric}
Let $u \in \mathcal{A}(n, \Lambda)$.  Suppose that $p \in \overline{\Omega} \cap B_{1/4}(0)$ and an $0<r\le 1/4$ such that $N_{\Omega}(p, r, u) \ge 2$. If either $p \in \partial \Omega$ or $B_r(p) \subset \Omega$ then, there is an $0< \epsilon= \epsilon(n, \Lambda)$ such that $u$ is not $(n-1, \epsilon, r, p)$-symmetric.
\end{lem}

\begin{proof}
Suppose that the statement is false.  Then, there exists a sequence of $u_i$, $p_{i}$, and $0<r_i$ such that $N_{T_{p_{i},r_i}\Omega_i}(0, 1, T_{p_{i}, r_i}u_i) \ge 2$ and $u_i$ is $(n-1, 2^{-i}, r, p)$-symmetric. By Lemma \ref{compactness}, we may extract a subsequence $T_{p_{j}, r_j}u_j$ which converges strongly in $W^{1, 2}(B_1(0))$ to a function $u_{\infty}.$ By Lemma \ref{N continuity}, $N_{\Omega_{\infty}}(0, 1, u_\infty) \ge 2$. But, $u_{\infty}$ is $(n-1, 0, 1, 0)$-symmetric and hence piece-wise linear.

If either $B_{1}(0) \subset \Omega_{\infty}$ or $0 \in \partial \Omega$ then $u_{\infty}$ being piecewise linear implies $N_{\Omega_{\infty}}(0, 1, u_\infty) =1$. This is a contradiction.  
\end{proof}

\subsection{Proof of Lemma \ref{e-reg containment}}
\begin{proof}
First, we prove (\ref{e:containment 1}).  Because for all integers $k$ and all $0< \epsilon$ the set $\mathcal{C}^k_{\epsilon}(u)$ is closed, we reduce to proving that there is an $0< \epsilon(n, \Lambda)$ such that $\mathcal{C}^{n-2}(u) \cap B_{1/8}(0) \subset \mathcal{C}^{n-2}_{\epsilon}(u)$. Suppose that this containment is false.  Then, there would exist a sequence of functions $u_i \in \mathcal{A}(n, \Lambda)$, points $p_{i} \in \mathcal{C}^{n-2}(u_i) \cap B_{\frac{1}{8}}(0)$ and scales $0< r_i \le 1$ such that $u_i$ is $(n-1, 2^{-i}, p_i, r_i)$-symmetric.  We rescale to the functions $T_{p_{i}, r_i}u_i$.  By Lemma \ref{strong convergence}, there exists a subsequence (also indexed by i) and a $(n-1)$-symmetric function $u_{\infty}$ such that $T_{p_i, r_i}u_i \rightarrow u_{\infty}$ strongly in $W^{1, 2}(B_1(0))$ and $C^{0}(B_1(0)).$  Note that for all $Q \in \partial \Omega_{\infty}$ and all $B_r(p) \subset \Omega_{\infty}$ $N_{\Omega_{\infty}}(Q, r, u_{\infty}) = N_{\Omega_{\infty}}(p, r, u_{\infty}) = 1$.

If $B_{\delta}(0) \in \Omega_{\infty}$ then for all sufficiently large $i \in \mathbb{N}$, $N_{\Omega_i}(p_{i}, \delta/2r_i ,u_i) \ge 2$. Letting $i \rightarrow \infty$, by Lemma \ref{N continuity}, $N_{\Omega_{\infty}}(0, \delta/2, u_{\infty})\ge 2$ which contradicts $u_{\infty}$ being piecewise linear.

If $0 \in \partial \Omega_{\infty}$, then for each $p_{i}$ we let $Q_{i} \in \partial \Omega_{i}$ be such that $|p_{i}- Q_{i}| = \text{dist}(p_{i}, \partial \Omega_{i})$. By Lemma \ref{force r small} there exists a $0< \delta(n, \Lambda)$ such that $N_{\Omega_i}(Q_{i}, 4|Q_{i} - p_{i}|, u_i) \ge 1+\delta$ for all $i \in \mathbb{N}$. For sufficiently large $i$ $4|Q_i-p_i| \le r_i$, and so by Lemma \ref{N monotonicity 1} for sufficiently large $i$ $N_{\Omega_i}(Q_{i}, r_i, u_i) \ge 1+\delta$.  Letting $i \rightarrow \infty$ we obtain that by Lemma \ref{N continuity} that $N_{\Omega_{\infty}}(0, 1, u_{\infty})> 1 + \delta$.  This is a contradiction. Thus, there exists an $0< \epsilon(n, \Lambda)$ such that $\mathcal{C}^{n-2}(u) \cap \Omega \cap B_{1/8}(0) \subset \mathcal{C}^{n-2}_{\epsilon}(u)$.

For prove \eqref{e:containment 2}, we note that if $Q \in \mathcal{C}^{n-2}(u) \cap \partial \Omega \setminus \text{sing}(\partial \Omega)$, then letting $r \rightarrow 0$, we may extract a subsequence such that $T_{Q, r_j}u \rightarrow u_{\infty}$ in the sense of Lemma \ref{compactness 1} and Lemma \ref{strong convergence}.  By the monotonicity of the Almgren frequency, Lemma \ref{N continuity}, and by considering $T_{Q, cr_j}u$ for any $0<c<1$ we see that $N_{\Omega_{\infty}}(0, r, u_{\infty}) \equiv \lim_{r \rightarrow 0^+}N_{\Omega}(Q, r, u)$.  Thus, Lemma \ref{N constant gives homogeneity} implies that $u_{\infty}$ is a homogeneous function which is harmonic in $\Omega_{\infty}$.  Since $Q \not \in \text{sing}(\partial \Omega)$, $\Omega_{\infty}$ is a half-space and we may extend $u_{\infty}$ to an entire, homogeneous harmonic function by reflection.  Since $Q \in \mathcal{C}^{n-2}(u)$ this polynomial must be a non-linear homogeneous harmonic polynomial and $\lim_{r \rightarrow 0^+}N_{\Omega}(Q, r, u) \ge 2$.  By Lemma \ref{N monotonicity 1}, then $N_{\Omega}(Q, r, u) \ge 2$ for all $0<r\le 1$ and Lemma \ref{boundary points not n-1 symmetric} gives the claim. 
\end{proof}

\section{Appendix A: H\"older Continuity}\label{S: Holder continuity}

In this section, we provide a proof of Lemma \ref{L: unif Holder bound 1}.  First, some standard results.

\begin{definition}
A bounded domain, $\Omega \subset \mathbb{R}^n$, is said to be of class $S$ if there exist numbers $0 < c_0 \le 1$ and $0< r_0$ such that for all $Q \in \partial \Omega$ and all $0< r\le r_0$
\begin{align*}
\mathcal{H}^{n}(B_r(Q) \cap \Omega^c) \ge c_0 \mathcal{H}^n(B_{r}(Q)).
\end{align*}
\end{definition}

\begin{lem}(Bounding the Supremum, \cite{Kenig94} Lemma 1.1.22)\label{sup bound}
Let $\Omega$ be a domain of class $S$. Let $Q \in \partial \Omega \cap B_1(0)$ and $0 < r \le \frac{1}{2}$.  Let $u$ be a function which is harmonic in $\Omega$ such that $u \in C(\overline{B_{2r}(Q) \cap \Omega})$ and $u \equiv 0$ on $B_{2r}(Q) \cap \partial \Omega$.  There exists a $c(n)$ such that for any $p \in B_r(Q) \cap \Omega$
\begin{align*}
\max_{B_r(Q) \cap \Omega}|u| \le c(n)\left( \fint_{B_{2r}(Q) \cap \Omega} u^2 dx \right)^{\frac{1}{2}}.
\end{align*}
\end{lem}

\begin{lem}(H\"older continuity up to the Boundary, \cite{Kenig94} Corollary 1.1.24)\label{holder boundary}
Let $\Omega$ be a domain of class $S$. Let $Q \in \partial \Omega \cap B_1(0)$ and $0 < r \le \frac{1}{2}$.  Let $u$ be a function which is harmonic in $\Omega$, $u \in C(\overline{B_{2r}(Q) \cap \Omega})$,  $u \equiv 0$ on $B_{2r}(Q) \cap \partial \Omega$, and $u \ge 0$.  There exists a $c(n)$ and an exponent $0 < \alpha (n) \le 1$ such that for any $p \in B_r(Q) \cap \Omega$
\begin{align*}
u(p) \le c(n) \left(\frac{|p - Q|}{r}\right)^{\alpha} \sup\{ u(y) : y \in B_{2r}(Q)\}.
\end{align*}
\end{lem}

\begin{thm}(Oscillation in the Interior, \cite{HeinonenKilpelainenMartio} Theorem 6.6)\label{T: interior osc}
Suppose that $u$ is harmonic in $\tilde \Omega$.  If $0 < r< R< \infty$ are such that $B_r(x) \subset B_R(x_0) \subset \tilde \Omega$, then
\begin{align*}
osc(u, B_r(x_0)) \le 2^{\alpha}\left(\frac{r}{R}\right)^{\alpha}osc(u, B_R(x_0)),
\end{align*}
where $\alpha = \alpha(n) \in (0, 1]$ only depends on $n$.
\end{thm}

\begin{thm}(H\"older Continuity in $B_2(0)$, \cite{HeinonenKilpelainenMartio} Theorem 6.44) \label{interior holder}
Suppose that $\Omega_1$ is of class $S$ with constant $c_0>0$ and $0< r_0 \le 1.$  Let $h \in C^0(\overline{\Omega_1})$ be a harmonic function in $\Omega_1$.  If there are constants $M \ge 0$ and $0 < \alpha \le 1$ such that 
\begin{align*}
|h(x) - h(y)| \le M |x -y|^{\alpha}
\end{align*}
for all $x, y \in \partial \Omega_1,$ then
\begin{align*}
|h(x) - h(y)| \le M_1 |x -y|^{\gamma} 
\end{align*}
for all $x, y \in \overline{\Omega_1}$.  Moreover, $\gamma = \gamma(n, \alpha, c_0)>0$ and one can choose $M_1 = 80 Mr_0^{-2}\max\{1, diam(\Omega_1)2\}$.
\end{thm}

\subsection{Proof of Lemma \ref{L: unif Holder bound 1}}
Let $u \in \mathcal{A}(n, \Lambda)$ with associated domain $\Omega \in \mathcal{D}(n)$, and let $Q_0 \in \partial \Omega \cap B_1(0)$ and $0 < r_0 \le \frac{1}{2}$.  First, we claim that  $(\fint_{B_2(0) \cap \Omega}(T_{Q_0, r_0}u)^2 dx)^{1/2} \le C(n, \Lambda).$ By the Poincare inequality
\begin{align*}
\fint_{B_2(0)} |T_{Q, r}u - Avg_{B_2(0)}(T_{Q, r}u)|^2 dx & \le C(n) |B_2(0)|^\frac{2}{n}\left(\fint_{B_2(0)} |\nabla T_{Q, r}u|^2dx\right)\\
& \le C(n)C(\Lambda, n).
\end{align*}
Furthermore, since $\mathcal{H}^{n}(\Omega^c \cap B_2(0)) \ge \frac{1}{10}\mathcal{H}^n(B_2(0)),$ we estimate $$|Avg_{B_2(0)}(T_{Q, r}u)|^2 \frac{1}{10}\mathcal{H}^n(B_2(0)) \le C(n)C(\Lambda, n).$$  Thus, $|\fint_{B_2(0)}(T_{Q, r}u)| \le C'(n) \sqrt{C(n , \Lambda)}.$

Next, we claim that there exists a $c(n)$ and an exponent $0 < \alpha (n) \le 1$ such that, for $Q \in T_{Q_0, r_0}\partial \Omega \cap B_1(0)$ and $p \in T_{Q_0, r_0}\Omega \cap B_{\frac{1}{2}}(Q)$
\begin{align*}
|T_{Q_0,r_0}u(p)| \le C(n, \Lambda) \left(\frac{|p - Q|}{r}\right)^{\alpha}.
\end{align*}

For any $T_{Q_0, r_0}u$ which changes sign in $T_{Q_0, r_0}\Omega \cap B_{1}(Q),$ we decompose $T_{Q_0, r_0}u = T_{Q_0, r_0}u^+ - T_{Q_0, r_0}u^-.$  Note that both $T_{Q_0, r_0}u^{\pm}$ are subharmonic.  Let $h_{\pm}$ be the harmonic extension of $T_{Q_0, r_0}u^{\pm}$ to $B_{1}(Q) \cap T_{Q_0, r_0}\Omega$.  Note that $B_{1}(Q) \cap T_{Q_0, r_0}\Omega$ is convex, and so is of class $S$.  Then, by Lemma \ref{holder boundary}, and the maximum principle
\begin{align*}
h_{\pm}(p) \le c(n) \left(\frac{|p - Q|}{r}\right)^{\alpha} \sup \{ h_{\pm}(y) : y \in \partial (B_{1}(Q) \cap T_{Q_0, r_0}\Omega) \}.
\end{align*}
By subharmonicity, $T_{Q_0, r_0}u^{\pm} \le h_{\pm}$, respectively.  By construction, $h_{\pm} = T_{Q_0, r_0}u^{\pm}$ on $\partial (B_{1}(Q) \cap T_{Q_0, r_0}\Omega)$ and by our first claim and Lemma \ref{sup bound}
\begin{align*}
\sup \{ h_{\pm}(y) : y \in \partial (B_{1}(Q) \cap T_{Q_0, r_0}\Omega)\} \le C(n, \Lambda).
\end{align*}
Note that this gives uniform control on the oscillation in $T_{Q_0, r_0}\Omega \cap B_2(0)$. This uniform control together with Theorem \ref{T: interior osc} implies that $T_{Q_0, r_0}u$ is locally H\"older on $\partial B_1(0) \cap T_{Q_0, r_0}\Omega$.

Now, we claim that for all $x, y \in \partial(T_{Q_0, r_0}\Omega \cap B_1(0))$
\begin{align*}
|T_{Q_0, r_0}u(x) - T_{Q_0, r_0}u(y)| \le C(n, \Lambda) |x -y|^{\alpha}.
\end{align*}
We argue by cases.  Suppose that $|x-y| < \max\{ \text{dist}(x, T_{Q_0, r_0}\partial \Omega), \text{dist}(y, T_{Q_0, r_0}\partial \Omega)\}$.  Then, there is a ball, $B_r(z) \subset T_{Q_0, r_0}\Omega$ with $|x-y| < r \le 2 |x-y|$ which contains both $x$ and $y$.  By Theorem \ref{T: interior osc} and the preceding paragraph, then we have the desired statement.

Suppose that $|x-y| \ge \max\{ \text{dist}(x, T_{Q_0, r_0}\partial \Omega), \text{dist}(y, T_{Q_0, r_0}\partial \Omega)\}$. Let $x_0, y_0 \in T_{Q_0, r_0}\partial \Omega$ be points such that $|x-x_0| = \text{dist}(x, \partial \Omega)$ and $|y-y_0| = \text{dist}(y, \partial \Omega)$.  Then
\begin{align*}
|T_{0, 1}u(x) - T_{0, 1}u(y)| & \le  |T_{0, 1}u(x) - T_{0, 1}u(x_0)| + |T_{0, 1}u(y) - T_{0, 1}u(y_0)|\\
& \le C(n, \Lambda) 2^{\alpha}|x- x_0|^{\alpha} + C(n, \Lambda) 2^{\alpha}|y- y_0|^{\alpha}\\
& \le C(n, \Lambda) 2^{\alpha +1} (\max\{ \text{dist}(x, T_{Q_0, r_0}\partial \Omega), \text{dist}(y, T_{Q_0, r_0}\partial \Omega)\})^{\alpha}\\
& \le C(n, \Lambda) |x -y|^{\alpha}.
\end{align*}
This proves the claim.  To obtain uniform interior H\"older continuity on the interior of $T_{Q_0, r_0}\Omega \cap B_1(0)$, we invoke Theorem \ref{interior holder} with $\Omega_1 = T_{Q_0, r_0}\Omega \cap B_1(0)$.

\bibliographystyle{alpha}
\bibliography{references}

\end{document}